\documentclass[12pt,a4paper]{article}
\usepackage{amsmath,amssymb,bm,ascmac,bbm,mathtools}
\usepackage[dvipdfmx, usenames]{color}
\usepackage[dvipdfmx]{graphicx}
\usepackage{xcolor}
\usepackage{here}
\usepackage{authblk}
\usepackage[hang,small,bf]{caption}
\usepackage[subrefformat=parens]{subcaption}
\usepackage{dcolumn}

\newcommand{\tr}{\text{tr}}

\newcommand{\mcal}{\mathcal}
\newcommand{\mbb}{\mathbb}
\newcommand{\mfrak}{\mathfrak}

\newcommand{\rank}{\text{rank}}

\newcommand{\wt}{\widetilde}

\newcommand{\mc}{\mathcal}

\newcommand{\wf}[1]{\widehat{\mfrak{#1}}}
\newcommand{\fp}{\text{FPdim}}
\newcommand{\vecG}{\text{Vec}}
\newcommand{\vect}{\text{Vect}_{\mbb C}}
\newcommand{\fib}{\text{Fib}}
\newcommand{\ising}{\text{Ising}}
\newcommand{\tc}{\text{ToricCode}}
\newcommand{\mods}{\text{mod }}
\newcommand{\Hom}{\underline{\text{Hom}}}

\begin{document}
\title{Classification of connected \'etale algebras in modular fusion categories up to rank five}
\author{Ken KIKUCHI}
\affil{Department of Physics, National Taiwan University, Taipei 10617, Taiwan}
\date{}
\maketitle

\begin{abstract}
We classify connected étale algebras in (possibly non-unitary) modular fusion categories $\mathcal B$'s with $\text{rank}(\mathcal B)\le5$. We also comment on Lagrangian algebra, anyon condensation, and physical applications. Concretely, we prove certain spontaneous $\mathcal B$-symmetry breaking and predict ground state degeneracies in massive renormalization group flows from non-unitary minimal models.
\end{abstract}

\makeatletter
\renewcommand{\theequation}
{\arabic{section}.\arabic{equation}}
\@addtoreset{equation}{section}
\makeatother

\section{Introduction}
Throughout the paper, $\mcal B$ denotes a modular fusion category (MFC) with braiding $c$ over the field $\mbb C$ of complex numbers (see \cite{EGNO15} for definitions). Its simple objects are denoted $b_i$'s with $i=1,\dots,\rank(\mcal B)$. (On the contrary, a fusion category (FC) over $\mbb C$ and its simple objects are denoted $\mcal C$ and $c_i$'s, respectively.) Thanks to the braiding, we can define an étale algebra $A\in\mcal B$. (The definitions are collected in section \ref{defs}.) The goal of this paper is to classify connected étale algebras in MFCs up to rank five. Our main results are summarized as the\newpage

\textbf{Theorem.} \textit{Connected étale algebras in modular fusion categories $\mcal B$'s up to rank five are given by}
\begin{table}[H]
\begin{center}
\makebox[1 \textwidth][c]{       
\resizebox{1.2 \textwidth}{!}{\begin{tabular}{c|c|c|c}
    Rank&$\mcal B$&Results&Completely anisotropic?\\\hline\hline
    1&$\vect$&e.g. \cite{KK23etale}&Yes\\\hline
    2&$\vecG_{\mbb Z/2\mbb Z}^{-1}$&\cite{KO01,KK23etale}&Yes\\
    &$\fib$&\cite{BD11,KK23etale}&Yes\\\hline
    3&$\vecG_{\mbb Z/3\mbb Z}^1$&\cite{EP09,KK23etale}&Yes\\
    &$\ising$&\cite{KO01,KL02,EM21,G23,KK23etale}&Yes\\
    &$psu(2)_5$&e.g. \cite{KK23etale}&Yes\\\hline
    4&$\vecG_{\mbb Z/2\mbb Z\times\mbb Z/2\mbb Z}^\alpha\simeq\begin{cases}\vecG_{\mbb Z/2\mbb Z}^{-1}\boxtimes\vecG_{\mbb Z/2\mbb Z}^{-1}\\\tc\end{cases}$&$\begin{cases}\text{Table }\ref{rank4vecGZ2Z2results}\\\text{Table }\ref{rank4toriccoderesults}\end{cases}$&$\begin{cases}\text{No }(\text{four with }d,h\text{ in }(\ref{Z2Z2etale}))/\text{ Yes (the other 16)}\\\text{No }(\text{four with }d,h\text{ in }(\ref{rank4ToricCoderesultalg}))/\text{ Yes (the other six)}\end{cases}$\\
    &$\vecG_{\mbb Z/4\mbb Z}^\alpha$&\cite{CEM23}, Table \ref{rank4Z4results}&Yes\\
    &$su(2)_3\simeq\vecG_{\mbb Z/2\mbb Z}^{-1}\boxtimes\fib$&\cite{KO01}, Table \ref{rank4su23results}&Yes\\
    &$\fib\boxtimes\fib$&\cite{BD11}, Table \ref{rank4fibfibresults}&No (four with $d,h$ in (\ref{fibfibetale}))/ Yes (the other 16)\\
    &$psu(2)_7$&Table \ref{rank4psu27results}&Yes\\\hline
    5&$su(5)_1$&\cite{G23}, Table \ref{rank5su51results}&Yes\\
    &$su(2)_4$&\cite{KO01}, Table \ref{rank5su24results}&No\\
    &$psu(2)_9$&Table \ref{rank5psu29results}&Yes\\
    &$psu(3)_4$&Table \ref{rank5psu34results}&Yes
\end{tabular}.}}
\end{center}
\caption{Connected étale algebras in MFC $\mcal B$ with $\rank(\mcal B)\le5$}\label{results}
\end{table}

\textbf{Remark.} For a reader's convenience, we also included known results up to rank three.\newline

\textbf{Remark.} Some MFCs are realized by Wess-Zumino-Witten (WZW) models. In those cases, we denote the MFCs sharing the same fusion ring by the realization, e.g., $su(2)_3$. Others are realized by subcategory of objects invariant under centers. We denote the MFCs by the realizations with $p$ in their head, e.g., $psu(2)_7$. Note that this is \textit{not} the same as `gauging' in physics language. For example, the $\mbb Z/2\mbb Z$ center symmetry in $su(2)_7$ is anomalous, and it cannot be gauged. (Mathematically, there is no $\mbb Z/2\mbb Z$ algebra.)

Especially in those cases, some classification results are known. (In this context, an MFC describing $\wf g_k$ WZW model is denoted $\mcal C(\mfrak g,k)$ and connected étale algebras in it is called quantum subgroups \cite{O00}. In another context, say \cite{CZW18}, connected étale algebras are also called condensable or normal algebras.) For example, connected étale algebras were classified in \cite{KO01} (for $\wf{su}(2)_k$), in \cite{EP09} (for $\wf{su}(3)_k$), in \cite{CEM23} (for $\wf{su}(4)_k$), and many more in \cite{G23}. When available, our results are consistent with them; $\mcal C(A_1,1),\mcal C(A_2,1),\mcal C(A_3,1),\mcal C(A_1,3),\mcal C(A_4,1)$ are known to be completely anisotropic, and $\mcal C(A_1,4)$ has two connected étale algebras $A\cong1,1\oplus X$ for $\mbb Z/2\mbb Z$ simple object $X\in\mcal C(A_1,4)$. The second gives the category of right $A$-modules corresponding to the $D_4$ Dynkin diagram \cite{KO01}, i.e., $\rank(\mcal C(A_1,4)_A)=4$.
\newline

\section{Classification}
\subsection{Definitions}\label{defs}
The monoidal products of the underlying fusion categories $\mcal C$'s are specified by fusion matrices $(N_i)_{jk}:={N_{i,j}}^k$ with $\mbb N$-coefficients
\begin{equation}
    c_i\otimes c_j\cong\bigoplus_{k=1}^{\rank(\mcal C)}{N_{i,j}}^kc_k.\label{fusioncoeff}
\end{equation}
Since entries of the fusion matrices are non-negative, we can apply the Perron-Frobenius theorem to get the largest eigenvalue $\fp_{\mcal C}(c_i)$ called the Frobenius-Perron dimension of $c_i$ (or $N_i$). Its squared sum defines Frobenius-Perron dimension of the fusion category
\begin{equation}
    \fp(\mcal C):=\sum_{i=1}^{\rank(\mcal C)}(\fp_{\mcal C}(c_i))^2.\label{FPdimB}
\end{equation}
In an MFC, one can also define quantum dimension $d_i$ of $c_i$ by the quantum (or categorical) trace
\[ d_i:=\tr(a_{c_i}), \]
where $a:id_{\mcal C}\cong(-)^{**}$ is a pivotal structure. They obey the same multiplication rules as simple objects
\begin{equation}
    d_id_j=\sum_{k=1}^{\rank(\mcal C)}{N_{i,j}}^kd_k.\label{quantumdim}
\end{equation}
Just as the Frobenius-Perron dimension, its squared sum defines categorical (or global) dimension
\begin{equation}
    D^2(\mcal C):=\sum_{i=1}^{\rank(\mcal C)}d_i^2.\label{categoricaldim}
\end{equation}
Note that for a given categorical dimension, there are two $D(\mcal C)$'s, one positive and one negative. A fusion category $\mcal C$ is called pseudo-unitary if $D^2(\mcal C)=\fp(\mcal C)$ and unitary if $\forall c_i\in\mcal C$, $d_i=\fp_{\mcal C}(c_i)$. These data characterize the underlying (spherical) fusion category.

An MFC has additional structure, braiding $c$. It is a natural isomorphism between two bifunctors
\begin{equation}
    c_{-,-}:-\otimes-\stackrel\sim\Rightarrow-\otimes-\label{braiding}
\end{equation}
subject to hexagon axioms. Here, two objects $b,b'$ (not necessarily simple) enter bifunctors with reversed orders:
\[ c_{b,b'}:b\otimes b'\stackrel\sim\Rightarrow b'\otimes b. \]
A fusion category with a braiding is called braided fusion category (BFC). It is specified by conformal dimensions $h_i$'s of $b_i$'s. For example, the double braiding of two simple objects $b_i,b_j$ is given by the formula
\begin{equation}
    c_{b_j,b_i}\cdot c_{b_i,b_j}\cong\sum_{k=1}^{\rank(\mcal B)}{N_{i,j}}^k\frac{e^{2\pi ih_k}}{e^{2\pi i(h_i+h_j)}}id_k,\label{doublebraiding}
\end{equation}
where $id_k$ is the identity morphism at $b_k$. Its quantum trace defines (unnormalized) $S$-matrix\footnote{Another modular matrix $T$ is also defined with conformal dimensions
\begin{equation}
    T_{i,j}:=e^{2\pi ih_i}\delta_{i,j}.\label{Tmatrix}
\end{equation}
The two modular matrices define additive central charge $c(\mcal B)$ mod 8 by
\begin{equation}
    (ST)^3=e^{2\pi ic(\mcal B)/8}S^2.\label{additivecentralcharge}
\end{equation}}
\begin{equation}
    \widetilde S_{i,j}:=\tr(c_{b_j,b_i}\cdot c_{b_i,b_j})=\sum_{k=1}^{\rank(\mcal B)}{N_{i,j}}^k\frac{e^{2\pi ih_k}}{e^{2\pi i(h_i+h_j)}}d_k.\label{Smatrix}
\end{equation}
A normalized $S$-matrix is defined by
\begin{equation}
    S_{i,j}:=\frac{\widetilde S_{i,j}}{D(\mcal B)}.\label{Snorm}
\end{equation}
An MFC is defined as a spherical BFC (called pre-modular fusion category) with non-degenerate $S$-matrix. It squares to charge conjugation
\[ S^2=C. \]
The charge conjugation matrix is defined by
\[ C_{i,j}=\delta_{i,j}\quad(b_i^*\cong b_j), \]
where $b_i^*$ is the dual of $b_i$. It obeys
\begin{equation}
    \widetilde S_{i,j^*}=\left(\widetilde S_{i,j}\right)^*,\label{Smatrixdual}
\end{equation}
where the RHS is the complex conjugate of $\widetilde S_{i,j}$. In particular, if $b_j$ is self-dual, $b_j^*\cong b_j$, then matrix elements should be real
\[ \forall b_i\in\mcal B,\quad\widetilde S_{i,j}\in\mbb R. \]
Using this condition, one can constrain conformal dimensions.

MFCs up to rank five (and those without multiplicity at rank six) have been classified in \cite{GK94,RSW07,BNRW15}. There are five rank four MFCs, and four rank five MFCs. We study all nine MFCs below.

Given the definitions on our ambient category, MFC, we next review definitions on algebras. An algebra in a fusion category $\mcal C$ is a triple $(A,\mu,u)$ of an object $A\in\mcal C$, multiplication morphism $\mu:A\otimes A\to A$, and unit morphism $u:1\to A$ obeying associativity and unit axioms. A category $\mcal C_A$ of right $A$-modules consists of pairs $(m,p)$ where $m\in\mcal C$ and $p:m\otimes A\to m$ subject to consistency conditions. (A category $_A\mcal C$ of left $A$-modules is defined analogously.) An algebra is called separable if $\mcal C_A$ is semisimple. An algebra $A\in\mcal B$ in a braided fusion category is called commutative if
\begin{equation}
    \mu\cdot c_{A,A}=\mu.\label{commutativealg}
\end{equation}
A commutative separable algebra is called étale. Any étale algebra $A$ decomposes to a direct sum of connected ones \cite{DMNO10} where $A\in\mcal C$ is called connected if $\dim_{\mbb C}\mcal C(1,A)=1$. A connected étale algebra $A\in\mcal B$ is called Lagrangian if $(\fp_{\mcal B}(A))^2=\fp(\mcal B)$. An example of connected étale algebra is the unit $A\cong1\in\mcal B$ giving $\mcal B_A\simeq\mcal B$. A BFC without nontrivial connected étale algebra is called completely anisotropic. The category of right $A$-modules has an important subcategory $\mcal B_A^0\subset\mcal B_A$. It consists of dyslectic (or local) modules \cite{P95} $(m,p)\in\mcal B_A$ obeying
\[ p\cdot c_{A,m}\cdot c_{m,A}=p. \]

The category $\mcal C_A$ of right $A$-modules is a left $\mcal C$-module category \cite{O01}. A left $\mcal C$-module category (or module category over $\mcal C$) is a quadruple $(\mcal M,\triangleright,m,l)$ of a category $\mcal M$, an action (or module product) bifunctor $\triangleright:\mcal C\times\mcal M\to\mcal M$, and natural isomorphisms $m_{-,-,-}:(-\otimes-)\triangleright-\cong-\triangleright(-\triangleright-)$ and $l:1\triangleright\mcal M\simeq\mcal M$ called module associativity constraint and unit constraint, respectively. They are subject to associativity and unit axioms. The actions form non-negative integer matrix representation (NIM-rep). This is because, in our setup, $\mc M$'s are fusion categories, and for any $b\in\mc B,m\in\mc M$, $b\triangleright m\in\mc M$ can be decomposed as a direct sum of simple objects with $\mbb N$-coefficients. The natural numbers assemble to $r\times r$ matrices where $r=\rank(\mc M)$. The NIM-rep is called $r$-dimensional.

\subsection{Method}
Here we explain our classification method. It is mainly based on \cite{KK23GSD}, however, as noticed in \cite{KK23etale}, modularity of our ambient category helps us. To explain this difference, we start from the method in \cite{KK23GSD}.

Let $\mcal B$ be a BFC and $A\in\mcal B$ a connected étale algebra. The Frobenius-Perron dimension of the category $\mcal B_A$ of right $A$-modules obeys \cite{KO01,ENO02,DMNO10}
\begin{equation}
    \fp_{\mcal B}(A)=\frac{\fp(\mcal B)}{\fp(\mcal B_A)}.\label{FPdimA}
\end{equation}
Since we have \cite{ENO02,EGNO15}
\begin{equation}
    \forall c\in\mcal C,\quad\fp_{\mcal C}(c)\ge1,\label{FPdimge1}
\end{equation}
we obtain an inequality
\begin{equation}
    1\le\fp_{\mcal B}(A)\le\fp(\mcal B).\label{FPdimAbound'}
\end{equation}
This led to a lemma \cite{KK23GSD}
\begin{equation}
    1\le\rank(\mcal B_A)\le\lfloor\fp(\mcal B)\rfloor.\label{rmax}
\end{equation}
Furthermore, since $A$ consists of simple objects of $\mcal B$, its general form is given by
\begin{equation}
    A\cong\bigoplus_{j=1}^{\rank(\mcal B)}n_jb_j,\label{generalA}
\end{equation}
where $n_j\in\mbb N$ counts the number of simple object $b_j$. Since the direct sum is defined as (co)limit, the object is equipped with product projections $p_j:A\to b_j$ and coproduct injections $\iota_j:b_j\to A$. Its Frobenius-Perron dimension is given by linear sum of those of simple objects:
\begin{equation}
    \fp_{\mcal B}(A)=\sum_{j=1}^{\rank(\mcal B)}n_j\fp_{\mcal B}(b_j).\label{FPdimAlinearsum}
\end{equation}
Based on these facts, in \cite{KK23GSD}, we first list candidate fusion categories $\mcal C$'s for $\mcal B_A$ obeying three necessary conditions: i) $1\le\rank(\mcal C)\le\lfloor\fp(\mcal B)\rfloor$, ii) $1\le\fp(\mcal C)\le\fp(\mcal B)$, and iii) $\fp_{\mcal B}(A)=\fp(\mcal B)/\fp(\mcal C)$. Then, we check which of the candidates satisfy axioms.

The method works nicely \cite{KK23GSD,KK23etale} for small $\fp(\mcal B)$ because fusion categories with small ranks are relatively well classified \cite{O02,O13,LPR20,VS22}. (See also \cite{L14,DZD16}.) The results are summarized in AnyonWiki \cite{anyonwiki}. On the other hand, the method heavily relies on lists of fusion categories, and works poorly for larger $\fp(\mcal B)$ where few fusion categories are known. However, if our ambient category $\mcal B$ is modular, we can do more.

The difference originates from the following fact \cite{P95,KO01,EGNO15}. Let $\mcal B$ be an MFC and $A\in\mcal B$ a connected étale algebra. Then, the category of dyslectic modules $\mcal B_A^0$ is modular obeying\footnote{The category of dyslectic modules further obeys
\[ e^{2\pi ic(\mcal B_A^0)/8}=e^{2\pi ic(\mcal B)/8}. \]}
\begin{equation}
    \fp(\mcal B_A^0)=\frac{\fp(\mcal B)}{\left(\fp_{\mcal B}(A)\right)^2}.\label{FPdimBA0}
\end{equation}
The same inequality (\ref{FPdimge1}) leads to
\begin{equation}
    1\le\left(\fp_{\mcal B}(A)\right)^2\le\fp(\mcal B).\label{fpdimAbound}
\end{equation}
Note the crucial difference of square. Because of the exponent, candidates of $A$ are more limited than the more general case above. Thus, we find only a few candidates without referring to lists of (modular) fusion categories. More concretely, in order to get connected algebra, we set an ansatz
\begin{equation}
    A\cong1\oplus\bigoplus_{b_j\not\cong1}n_jb_j\label{ansatz}
\end{equation}
with $n_j\in\mbb N$. It has
\[ \fp_{\mcal B}(A)=1+\sum_{b_j\not\cong1}n_j\fp_{\mcal B}(b_j). \]
Plugging this into (\ref{fpdimAbound}), we usually find only a few sets of $n_j$'s. We check which set with (\ref{ansatz}) satisfies the axioms. (In this way, we could prove \cite{KK23etale} all $psu(2)_5$'s are completely anisotropic.)

In checking axioms, we have to make sure three points; $A$ be connected, separable, and commutative. The connectedness is automatic from our ansatz (\ref{ansatz}). The separability is also automatic if $\mcal B_A$ is a fusion category. Thus, to check separability, we check quantum dimensions of right $A$-modules $m\in\mc B_A$ \cite{KO01}
\begin{equation}
    d_{\mc B_A}(m)=\frac{d_{\mc B}(m)}{d_{\mc B}(A)}.\label{dBAm}
\end{equation}
where $d_{\mcal C}(X)$ is the quantum dimension of $X$ viewed as an object of $\mcal C$. What is nontrivial is commutativity of $A$. In order to check the axiom (\ref{commutativealg}), we in principle have to compute braiding $c_{A,A}$. Namely, we have to solve hexagon equations. This task is in general cumbersome. So as to reduce our task, we find it useful to consider a necessary condition
\begin{equation}
    \mu\cdot c_{A,A}\cdot c_{A,A}=\mu\label{commutativenecessary}
\end{equation}
for an algebra to be commutative. Since the condition is written in terms of double braiding, we can employ the formula (\ref{doublebraiding}). Concretely, the double braiding is given by
\begin{equation}
    c_{A,A}\cdot c_{A,A}\cong\sum_{i,j=1}^{\rank(\mcal B)}n_in_j(\iota_i\otimes\iota_j)\cdot c_{b_j,b_i}\cdot c_{b_i,b_j}\cdot(p_i\otimes p_j).\label{doublecAA}
\end{equation}
If a candidate fails to satisfy the necessary condition, we can safely discard the candidate. In this way, we reduce the number of candidates. Below, we demonstrate our method in examples with $\rank(\mcal B)\le5$. (Those with $\rank(\mcal B)\le3$ have been studied in \cite{KK23etale}.)

\subsection{Rank four}
\subsubsection{$\mcal B\simeq\vecG_{\mbb Z/2\mbb Z\times\mbb Z/2\mbb Z}^\alpha$}
The MFCs have four simple objects $\{1,X_{h_X},Y_{h_Y},Z_{h_Z}\}$ obeying monoidal products
\begin{table}[H]
\begin{center}
\begin{tabular}{c|c|c|c|c}
    $\otimes$&1&$X$&$Y$&$Z$\\\hline
    1&1&$X$&$Y$&$Z$\\\hline
    $X$&&$1$&$Z$&$Y$\\\hline
    $Y$&&&1&$X$\\\hline
    $Z$&&&&$1$
\end{tabular}.
\end{center}
\end{table}
\hspace{-17pt}While the monoidal products are the same, there are two classes of MFCs; one is given by a Deligne tensor product \cite{D90} $\mcal B\simeq\vecG_{\mbb Z/2\mbb Z}^{-1}\boxtimes\vecG_{\mbb Z/2\mbb Z}^{-1}$, and another is not a product $\mcal B\simeq\tc$. (Following \cite{RSW07}, we call the latter the Toric Code MFCs \cite{K96,K97}.) Their difference appears in conformal dimensions. Thus, up to there, we study them simultaneously.

Because of the monoidal products, they have
\[ \fp_{\mcal B}(1)=\fp_{\mcal B}(X)=\fp_{\mcal B}(Y)=\fp_{\mcal B}(Z)=1, \]
and
\[ \fp(\mcal B)=4. \]
Their quantum dimensions $d_j$'s are solutions of
\[ d_X^2=1,\quad d_Xd_Y=d_Z,\quad d_Xd_Z=d_Y,\quad d_Y^2=1,\quad d_Yd_Z=d_X,\quad d_Z^2=1. \]
There are four solutions
\[ (d_X,d_Y,d_Z)=(1,1,1),(1,-1,-1),(-1,1,-1),(-1,-1,1). \]
All solutions have the same categorical dimension
\[ D^2(\mcal B)=4. \]
Only the first solution is unitary.

In order to count the number of MFCs with the fusion ring, let us compute conformal dimensions. Here, we have to treat two classes $\mcal B\simeq\vecG_{\mbb Z/2\mbb Z}^{-1}\boxtimes\vecG_{\mbb Z/2\mbb Z}^{-1},\tc$ separately. We start from the product MFCs.

\paragraph{$\mcal B\simeq\vecG_{\mbb Z/2\mbb Z}^{-1}\boxtimes\vecG_{\mbb Z/2\mbb Z}^{-1}$.} For our $\mcal B$ to be modular, we need both factors $\vecG_{\mbb Z/2\mbb Z}^\alpha$'s be modular. Since each factor has $S$-matrix
\[ \wt S=\begin{pmatrix}1&d_b\\d_b&e^{-4\pi ih_b}\end{pmatrix} \]
for the $\mbb Z/2\mbb Z$ object $b$, we need
\begin{equation}
    h_b=\frac14\quad(\mods\frac12).\label{confdimZ2}
\end{equation}
Denote the $\mbb Z/2\mbb Z$ objects of two factors $X$ and $Y$. Modular $\mcal B$ is given by
\[ (h_X,h_Y)=(\frac14,\frac14),(\frac14,\frac34),(\frac34,\frac34)\quad(\mods1). \]
(This implies both factors should be anomalous $\alpha=-1=\alpha'$. We have already specified this fact in the name of the MFCs.) Since conformal dimension of a Deligne tensor product $Z\cong X\boxtimes Y$ is given by $h_Z=h_X+h_Y$ mod 1, we have
\[ h_Z=\frac12,0,\frac12\quad(\mods1), \]
respectively. (This means the diagonal $\mbb Z/2\mbb Z$ generated by $Z$ is anomaly-free.) How many MFCs are there? Note that we should not distinguish two MFCs if they are simply related by exchange of two factors (or changing names of simple objects). In order to avoid the double-counting, we perform case analysis.

\begin{itemize}
\item $(d_X,d_Y)=(1,1)$. This gives unitary MFCs. Different MFCs are given by conformal dimensions
\[ (h_X,h_Y)=(\frac14,\frac14),(\frac14,\frac34),(\frac34,\frac34)\quad(\mods1). \]
Including two signs of categorical dimensions, we have six MFCs, consistent with \cite{RSW07}.
\item $(d_X,d_Y)=(1,-1)$. In this case, different MFCs are given by conformal dimensions
\[ (h_X,h_Y)=(\frac14,\frac14),(\frac14,\frac34),(\frac34,\frac14),(\frac34,\frac34)\quad(\mods1). \]
Taking two categorical dimensions into account, we have eight MFCs.
\item $(d_X,d_Y)=(-1,-1)$. Different MFCs are given by conformal dimensions
\[ (h_X,h_Y)=(\frac14,\frac14),(\frac14,\frac34),(\frac34,\frac34)\quad(\mods1). \]
With two signs of categorical dimensions, we have six MFCs.
\end{itemize}
Therefore, there are
\[ 6+8+6=20 \]
product MFCs in total. Among them, those six with quantum dimensions $+1$ are unitary. Note that quantum dimension of $Z$ would also matter in classifying connected étale algebra. (We will see only $Z$ can enter the algebra.) Its quantum dimension is given by the product $d_Z=d_Xd_Y$. Therefore, the $\mbb Z/2\mbb Z$ object $Z$ is a boson iff $(d_Z,h_Z)=(1,0),(-1,\frac12)$ (mod 1 for conformal dimensions) \cite{KK23etale}, or more concretely
\begin{equation}
    c_{Z,Z}\cong id_1\iff(d_X,d_Y,h_X,h_Y)=(1,1,\frac14,\frac34),(-1,-1,\frac14,\frac34),(1,-1,\frac14,\frac14),(1,-1,\frac34,\frac34)\quad(\mods1\text{ for }h).\label{cZZtrivial}
\end{equation}

\paragraph{$\mcal B\simeq\tc$.} In these MFCs, the simple objects have conformal dimensions $0$ or $\frac12$ mod 1. Studying necessary conditions, regardless of quantum dimensions, one finds four sets
\[ (h_X,h_Y,h_Z)=(\frac12,0,0),(0,\frac12,0),(0,0,\frac12),(\frac12,\frac12,\frac12)\quad(\mods1). \]
In order to count the number of MFCs, we again perform case analysis.
\begin{itemize}
    \item $(d_X,d_Y,d_Z)=(1,1,1)$. This gives unitary MFCs. Different MFCs are given by conformal dimensions
    \[ (h_X,h_Y,h_Z)=(\frac12,0,0),(\frac12,\frac12,\frac12)\quad(\mods1). \]
    Taking two signs of categorical dimensions into account, we find four unitary Toric Code MFCs, consistent with \cite{RSW07}.
    \item $(d_X,d_Y,d_Z)=(1,-1,-1)$. Different MFCs are given by conformal dimensions
    \[ (h_X,h_Y,h_Z)=(\frac12,0,0),(0,\frac12,0),(\frac12,\frac12,\frac12)\quad(\mods1). \]
    With two signs of categorical dimensions, we find six MFCs.
\end{itemize}
Therefore, there are
\[ 4+6=10 \]
Toric Code MFCs in total, among which those four with positive quantum dimensions give unitary MFCs.

Including both $\mcal B\simeq\vecG_{\mbb Z/2\mbb Z}^{-1}\boxtimes\vecG_{\mbb Z/2\mbb Z}^{-1}$ and $\mcal B\simeq\tc$, we find
\[ 20+10=30 \]
MFCs with $\mbb Z/2\mbb Z\times\mbb Z/2\mbb Z$ fusion ring. Among them, six from the product MFCs and four from the Toric Code MFCs are unitary, ten in total, consistent with \cite{RSW07}.

Having computed conformal dimensions, we can classify connected étale algebras. The most general ansatz
\[ A\cong1\oplus n_XX\oplus n_YY\oplus n_ZZ \]
with $n_j\in\mbb N$ has
\[ \fp_{\mcal B}(A)=1+n_X+n_Y+n_Z. \]
For this to obey (\ref{fpdimAbound}), the only possible sets are given by
\[ (n_X,n_Y,n_Z)=(0,0,0),(1,0,0),(0,1,0),(0,0,1). \]
The first candidate is nothing but the trivial connected étale algebra $A\cong1$ giving $\mcal B_A^0\simeq\mcal B_A\simeq\mcal B$. Whether the other three can be connected étale or not depends on quantum and conformal dimensions.

For $\mcal B\simeq\vecG_{\mbb Z/2\mbb Z}^{-1}\boxtimes\vecG_{\mbb Z/2\mbb Z}^{-1}$, $X,Y$ have nontrivial conformal dimensions, and candidates with them fail to be commutative because they do not satisfy the necessary condition (\ref{commutativenecessary}). On the other hand, the last candidate with $Z$ can be commutative iff $d,h$ are given as (\ref{cZZtrivial}). What remains is to check the separability. Let us study this point by identifying $\mc B_A$. Its $\fp_{\mcal B}(A)=2$ demands
\[ \fp(\mcal B_A^0)=1,\quad\fp(\mcal B_A)=2. \]
The only possibility is
\[ \mcal B_A^0\simeq\vect,\quad\mcal B_A\simeq\vecG_{\mbb Z/2\mbb Z}^\alpha. \]
Since the nontrivial simple object of $\vecG_{\mbb Z/2\mbb Z}^\alpha$ have quantum dimensions $\pm1$, a commutative algebra (\ref{cZZtrivial}) is separable iff $d_Z=1$. In other words, we get connected étale algebra
\begin{equation}
    A\cong1\oplus Z\quad(d_X,d_Y,h_X,h_Y)=(1,1,\frac14,\frac34),(-1,-1,\frac14,\frac34)\quad(\mods1\text{ for }h).\label{Z2Z2etale}
\end{equation}

For $\mcal B\simeq\tc$, the nontrivial simple objects with $(d,h)=(1,0),(-1,\frac12)$ give commutative algebras. They are separable when $d=1$ as above. Therefore, we find four connected étale algebras:
\begin{equation}
    \hspace{-40pt}A\cong\begin{cases}1&(\text{all MFCs}),\\1\oplus X&(d_X,d_Y,d_Z,h_X,h_Y,h_Z)=(1,-1,-1,0,\frac12,0),\\1\oplus Y&(d_X,d_Y,d_Z,h_X,h_Y,h_Z)=(1,1,1,\frac12,0,0),\\1\oplus Z&(d_X,d_Y,d_Z,h_X,h_Y,h_Z)=(1,1,1,\frac12,0,0).\end{cases}\quad(\mods1\text{ for }h)\label{rank4ToricCoderesultalg}
\end{equation}

To summarize, we found
\begin{table}[H]
\begin{center}
\begin{tabular}{c|c|c|c}
    Connected étale algebra $A$&$\mcal B_A$&$\rank(\mcal B_A)$&Lagrangian?\\\hline
    1&$\mcal B$&4&No\\
    $1\oplus Z$ for (\ref{Z2Z2etale})&$\vecG_{\mbb Z/2\mbb Z}^\alpha$&2&Yes
\end{tabular},
\end{center}
\caption{Connected étale algebras in rank four MFC $\mcal B\simeq\vecG_{\mbb Z/2\mbb Z}^{-1}\boxtimes\vecG_{\mbb Z/2\mbb Z}^{-1}$}\label{rank4vecGZ2Z2results}
\end{table}
\hspace{-17pt}and
\begin{table}[H]
\begin{center}
\begin{tabular}{c|c|c|c}
    Connected étale algebra $A$&$\mcal B_A$&$\rank(\mcal B_A)$&Lagrangian?\\\hline
    1&$\mcal B$&4&No\\
    $1\oplus X$ for (\ref{rank4ToricCoderesultalg})&$\vecG_{\mbb Z/2\mbb Z}^\alpha$&2&Yes\\
    $1\oplus Y$ for (\ref{rank4ToricCoderesultalg})&$\vecG_{\mbb Z/2\mbb Z}^\alpha$&2&Yes\\
    $1\oplus Z$ for (\ref{rank4ToricCoderesultalg})&$\vecG_{\mbb Z/2\mbb Z}^\alpha$&2&Yes
\end{tabular}.
\end{center}
\caption{Connected étale algebras in rank four MFC $\mcal B\simeq\tc$}\label{rank4toriccoderesults}
\end{table}
\hspace{-17pt}In particular, 16 product MFCs with
\begin{align*}
    (d_X,d_Y,h_X,h_Y)=&(1,1,\frac14,\frac14),(1,1,\frac34,\frac34),(1,-1,\frac14,\frac14),(1,-1,\frac14,\frac34),\\
    &(1,-1,\frac34,\frac14),(1,-1,\frac34,\frac34),(-1,-1,\frac14,\frac14),(-1,-1,\frac34,\frac34)\quad(\mods1\text{ for }h)
\end{align*}
and six Toric Code MFCs with
\[ (d_X,d_Y,d_Z,h_X,h_Y,h_Z)=(1,1,1,\frac12,\frac12,\frac12),(1,-1,-1,\frac12,0,0),(1,-1,-1,\frac12,\frac12,\frac12)\quad(\mods1\text{ for }h) \]
are completely anisotropic.
\newline

\subsubsection{$\mcal B\simeq\vecG_{\mbb Z/4\mbb Z}^\alpha$}
The MFCs have four simple objects $\{1,X_{h_X},Y_{h_Y},Z_{h_Z}\}$ obeying monoidal products
\begin{table}[H]
\begin{center}
\begin{tabular}{c|c|c|c|c}
    $\otimes$&1&$X$&$Y$&$Z$\\\hline
    1&1&$X$&$Y$&$Z$\\\hline
    $X$&&$Y$&$Z$&1\\\hline
    $Y$&&&1&$X$\\\hline
    $Z$&&&&$Y$
\end{tabular}.
\end{center}
\end{table}
\hspace{-17pt}Thus, they have
\[ \fp_{\mcal B}(1)=\fp_{\mcal B}(X)=\fp_{\mcal B}(Y)=\fp_{\mcal B}(Z)=1, \]
and
\[ \fp(\mcal B)=4. \]
Their quantum dimensions $d_j$'s are solutions of
\[ d_X^2=d_Y,\quad d_Xd_Y=d_Z,\quad d_Xd_Z=1,\quad d_Y^2=1,\quad d_Yd_Z=d_X,\quad d_Z^2=d_Y. \]
There are two solutions
\[ (d_X,d_Y,d_Z)=(1,1,1),(-1,1,-1). \]
The first solution gives unitary MFCs and the second gives non-unitary ones. They both have the same categorical dimension
\[ D^2(\mcal B)=4. \]

In order to count the number of MFCs, we compute their conformal dimensions. Studying some necessary conditions such as non-degeneracy of $S$-matrix or $S^2=C$ where $C$ is the charge conjugation matrix, regardless of quantum dimensions, we find four conformal dimensions
\[ (h_X,h_Y,h_Z)=(\frac18,\frac12,\frac18),(\frac38,\frac12,\frac38),(\frac58,\frac12,\frac58),(\frac78,\frac12,\frac78)\quad(\mods1). \]
(This implies the $\mbb Z/4\mbb Z$ is anomalous.) Therefore, there are
\[ 2(\text{quantum dimensions})\times4(\text{conformal dimensions})\times2(\text{categorical dimensions})=16 \]
MFCs, among which those eight with positive quantum dimensions give unitary MFCs, consistent with \cite{RSW07}. We classify connected étale algebras in all 16 MFCs simultaneously.

The most general ansatz for a connected algebra
\[ A\cong1\oplus n_XX\oplus n_YY\oplus n_ZZ \]
with $n_j\in\mbb N$ has
\[ \fp_{\mcal B}(A)=1+n_X+n_Y+n_Z. \]
For this to obey (\ref{fpdimAbound}), the natural numbers are restricted to one of four
\[ (n_X,n_Y,n_Z)=(0,0,0),(1,0,0),(0,1,0),(0,0,1). \]
The first is nothing but the trivial connected étale algebra giving $\mcal B_A^0\simeq\mcal B_A\simeq\mcal B$. The second and the fourth cannot be commutative because $X,Z$ have nontrivial conformal dimensions. The one with the $\mbb Z/2\mbb Z$ object $Y$ also fails to be commutative because it cannot have $(d_Y,h_Y)=(1,0),(-1,\frac12)$.

To summarize, we arrive
\begin{table}[H]
\begin{center}
\begin{tabular}{c|c|c|c}
    Connected étale algebra $A$&$\mcal B_A$&$\rank(\mcal B_A)$&Lagrangian?\\\hline
    1&$\mcal B$&4&No
\end{tabular}.
\end{center}
\caption{Connected étale algebras in rank four MFC $\mcal B\simeq\vecG_{\mbb Z/4\mbb Z}^\alpha$}\label{rank4Z4results}
\end{table}
\hspace{-17pt}Namely, all 16 MFCs $\vecG_{\mbb Z/4\mbb Z}^\alpha$'s are completely anisotropic.
\newline

\subsubsection{$\mcal B\simeq su(2)_3\simeq\vecG_{\mbb Z/2\mbb Z}^{-1}\boxtimes\fib$}
The MFCs have four simple objects $\{1,X_{h_X},Y_{h_Y},Z_{h_Z}\}$ obeying monoidal products
\begin{table}[H]
\begin{center}
\begin{tabular}{c|c|c|c|c}
    $\otimes$&1&$X$&$Y$&$Z$\\\hline
    1&1&$X$&$Y$&$Z$\\\hline
    $X$&&$1$&$Z$&$Y$\\\hline
    $Y$&&&$1\oplus Y$&$X\oplus Z$\\\hline
    $Z$&&&&$1\oplus Y$
\end{tabular}.
\end{center}
\end{table}
\hspace{-17pt}Thus, they have
\[ \fp_{\mcal B}(1)=1=\fp_{\mcal B}(X),\quad\fp_{\mcal B}(Y)=\zeta:=\frac{1+\sqrt5}2=\fp_{\mcal B}(Z), \]
and
\[ \fp(\mcal B)=5+\sqrt5\approx7.2. \]
The quantum dimensions are solutions of $d_X^2=1,d_Xd_Y=d_Z,d_Xd_Z=d_Y,d_Y^2=1+d_Y,d_Yd_Z=d_X+d_Z,d_Z^2=1+d_Y$. There are four solutions
\[ (d_X,d_Y,d_Z)=(1,\zeta,\zeta),(-1,\zeta,-\zeta),(1,-\zeta^{-1},-\zeta^{-1}),(-1,-\zeta^{-1},\zeta^{-1}). \]
(Only the first solution gives unitary MFCs.) Accordingly, there are two categorical dimensions
\[ D^2(\mcal B)=5-\sqrt5(\approx2.8),\quad5+\sqrt5. \]
For $su(2)_3\simeq\vecG_{\mbb Z/2\mbb Z}^\alpha\boxtimes\fib$ to be modular, both factors should be modular. Thus, conformal dimensions $h_X,h_Y$ (and $h_Z=h_X+h_Y$ mod 1) can only take the following eight values
\[ (h_X,h_Y,h_Z)=\begin{cases}(\frac14,\frac25,\frac{13}{20}),(\frac14,\frac35,\frac{17}{20}),(\frac34,\frac25,\frac3{20}),(\frac34,\frac35,\frac7{20})&(d_Y=\zeta),\\(\frac14,\frac15,\frac9{20}),(\frac14,\frac45,\frac1{20}),(\frac34,\frac15,\frac{19}{20}),(\frac34,\frac45,\frac{11}{20})&(d_Y=-\zeta^{-1}).\end{cases}\quad(\mods1) \]
Therefore, there are
\[ 2^2(\text{quantum dimensions})\times4(\text{conformal dimensions})\times2(\text{categorical dimensions})=32 \]
MFCs. Among them, those eight with $(d_X,d_Y)=(+1,\zeta)$ (and $d_Z=d_Xd_Y=\zeta$) are unitary MFCs, consistent with \cite{RSW07}. We classify connected étale algebras in all 32 MFCs simultaneously.

Set an ansatz
\[ A\cong1\oplus n_XX\oplus n_YY\oplus n_ZZ \]
with $n_X,n_Y,n_Z\in\mbb N$. It has
\[ \fp_{\mcal B}(A)=1+n_X+(n_Y+n_Z)\zeta. \]
For this to obey (\ref{fpdimAbound}), the coefficients can take only four values
\[ (n_X,n_Y,n_Z)=(0,0,0),(1,0,0),(0,1,0),(0,0,1). \]
The first possibility is nothing but the trivial connected étale algebra $A\cong1$ giving $\mcal B_A^0\simeq\mcal B_A\simeq\mcal B$. The other three candidates fail to satisfy the necessary condition (\ref{commutativenecessary}) because nontrivial simple objects have nontrivial conformal dimensions. Thus, the three candidates do not give connected étale algebra.

We conclude
\begin{table}[H]
\begin{center}
\begin{tabular}{c|c|c|c}
    Connected étale algebra $A$&$\mcal B_A$&$\rank(\mcal B_A)$&Lagrangian?\\\hline
    1&$\mcal B$&4&No
\end{tabular}.
\end{center}
\caption{Connected étale algebras in rank four MFC $\mcal B\simeq su(2)_3$}\label{rank4su23results}
\end{table}
\hspace{-17pt}Namely, all 32 $\mcal B\simeq su(2)_3$'s are completely anisotropic.

\subsubsection{$\mcal B\simeq\fib\boxtimes\fib$}
The MFCs have four simple objects $\{1,X_{h_X},Y_{h_Y},Z_{h_Z}\}$ obeying monoidal products
\begin{table}[H]
\begin{center}
\begin{tabular}{c|c|c|c|c}
    $\otimes$&1&$X$&$Y$&$Z$\\\hline
    1&1&$X$&$Y$&$Z$\\\hline
    $X$&&$1\oplus X$&$Z$&$Y\oplus Z$\\\hline
    $Y$&&&$1\oplus Y$&$X\oplus Z$\\\hline
    $Z$&&&&$1\oplus X\oplus Y\oplus Z$
\end{tabular}.
\end{center}
\end{table}
\hspace{-17pt}Thus, they have
\[ \fp_{\mcal B}(1)=1,\quad\fp_{\mcal B}(X)=\zeta=\fp_{\mcal B}(Y),\quad\fp_{\mcal B}(Z)=\frac{3+\sqrt5}2, \]
and
\[ \fp(\mcal B)=\frac{15+5\sqrt5}2\approx13.1. \]

The Fibonacci objects $X,Y$ can have quantum dimensions
\[ d=\zeta,-\zeta^{-1}. \]
Accordingly, they have categorical dimensions
\[ D^2(\mcal B)=\frac{15-5\sqrt5}2(\approx1.9),\quad5,\quad\frac{15+5\sqrt5}2. \]
Depending on the quantum dimensions, they can have conformal dimensions
\[ h=\begin{cases}\frac25,\frac35&(d=\zeta),\\\frac15,\frac45&(d=-\zeta^{-1}).\end{cases}\quad(\mods1) \]
(The last nontrivial simple object $Z\cong X\boxtimes Y$ has $d_Z=d_Xd_Y$ and $h_Z=h_X+h_Y$ mod 1.) How many MFCs are there? As in $\mcal B\simeq\vecG_{\mbb Z/2\mbb Z}^{-1}\boxtimes\vecG_{\mbb Z/2\mbb Z}^{-1}$, we should not distinguish two MFCs if they are related by exchange of two factors. In order to avoid the double-counting, we perform case analysis.

\paragraph{$(d_X,d_Y)=(\zeta,\zeta)$.} All distinct MFCs with these quantum dimensions are labeled by conformal dimensions
\[ (h_X,h_Y)=(\frac25,\frac25),(\frac25,\frac35),(\frac35,\frac35)\quad(\mods1). \]
In particular, we should not count $(h_X,h_Y)=(\frac35,\frac25)$ because the MFC with the conformal dimensions is isomorphic to the one with $(h_X,h_Y)=(\frac25,\frac35)$. Therefore, there are
\[ 3(\text{conformal dimensions})\times2(\text{categorical dimensions})=6 \]
MFCs with this quantum dimension. Since these are unitary, the number is consistent with \cite{RSW07}.

\paragraph{$(d_X,d_Y)=(\zeta,-\zeta^{-1})$.} All distinct MFCs with these quantum dimensions are labeled by conformal dimensions
\[\ (h_X,h_Y)=(\frac25,\frac15),(\frac25,\frac45),(\frac35,\frac15),(\frac35,\frac45)\quad(\mods1). \]
Thus, there are
\[ 4(\text{conformal dimensions})\times2(\text{categorical dimensions})=8 \]
MFCs.

\paragraph{$(d_X,d_Y)=(-\zeta^{-1},-\zeta^{-1})$.} All distinct MFCs with these quantum dimensions are labeled by conformal dimensions
\[\ (h_X,h_Y)=(\frac15,\frac15),(\frac15,\frac45),(\frac45,\frac45)\quad(\mods1). \]
Thus, there are
\[ 3(\text{conformal dimensions})\times2(\text{categorical dimensions})=6 \]
MFCs.

In total, there are
\[ 6+8+6=20 \]
MFCs $\mcal B\simeq\fib\boxtimes\fib$, among which six are unitary. We classify connected étale algebras in all 20 MFCs simultaneously.

We work with an ansatz
\[ A\cong1\oplus n_XX\oplus n_YY\oplus n_ZZ \]
with $n_X,n_Y,n_Z\in\mbb N$. It has
\[ \fp_{\mcal B}(A)=1+(n_X+n_Y)\zeta+n_Z\frac{3+\sqrt5}2. \]
For this to obey (\ref{fpdimAbound}), there are only four possibilities
\[ (n_X,n_Y,n_Z)=(0,0,0),(1,0,0),(0,1,0),(0,0,1). \]
The first candidate is nothing but the trivial connected étale algebra $A\cong1$ giving $\mcal B_A^0\simeq\mcal B_A\simeq\mcal B$. The second and third candidates with $X$ or $Y$ do not pass the necessary condition (\ref{commutativenecessary}) because $X$ and $Y$ have nontrivial conformal dimensions. The last candidate with $Z$ can pass the necessary condition because $Z$ can have trivial conformal dimension. Concretely, $Z$ has trivial conformal dimension when $d_X=d_Y$ and $h_X\neq h_Y$ mod 1. Indeed, it is known \cite{BD11} that the candidate is a commutative algebra when the two factor of $\fib$'s are reverse (or opposite) to each other. It further turns out separable, hence étale. To check this point, we identify $\mc B_A$.

The candidate $A\cong1\oplus Z$ has
\[ \fp_{\mcal B}(A)=\frac{5+\sqrt5}2, \]
and demands
\[ \fp(\mcal B_A^0)=1,\quad\fp(\mc B_A)=\frac{5+\sqrt5}2. \]
The only possible category of dyslectic (right) $A$-modules is $\mcal B_A^0\simeq\vect$. (This identification also matches central charges.) The category $\mc B_A$ of right $A$-modules is identified as
\[ \mc B_A\simeq\fib. \]
To show this, we search for NIM-reps. Indeed, we find a two-dimensional NIM-rep
\[ n_1=1_2,\quad n_X=\begin{pmatrix}0&1\\1&1\end{pmatrix}=n_Y,\quad n_Z=\begin{pmatrix}1&1\\1&2\end{pmatrix}. \]
Denoting a basis of $\mc M$ by $\{m_1,m_2\}$, we obtain identifications
\[ F(m_1)\cong1\oplus Z,\quad F(m_2)\cong X\oplus Y\oplus Z. \]
They have quantum dimensions (\ref{dBAm})
\[ d_{\mc B_A}(F(m_1))=1,\quad d_{\mc B_A}(F(m_2))=d_X. \]
Together with its Frobenius-Perron dimension, we get the identification. Since $\mc B_A\simeq\fib$ is semisimple, $A$ is separable and étale. We obtain
\begin{equation}
    A\cong1\oplus Z\quad(d_X,d_Y,h_X,h_Y)=(\zeta,\zeta,\frac25,\frac35),(-\zeta^{-1},-\zeta^{-1},\frac15,\frac45).\quad(\mods1\text{ for }h)\label{fibfibetale}
\end{equation}

To summarize, we conclude
\begin{table}[H]
\begin{center}
\begin{tabular}{c|c|c|c}
    Connected étale algebra $A$&$\mcal B_A$&$\rank(\mcal B_A)$&Lagrangian?\\\hline
    1&$\mcal B$&4&No\\
    $1\oplus Z$ for (\ref{fibfibetale})&$\fib$&2&Yes
\end{tabular}.
\end{center}
\caption{Connected étale algebras in rank four MFC $\mcal B\simeq\fib\boxtimes\fib$}\label{rank4fibfibresults}
\end{table}
\hspace{-17pt}That is, those four in (\ref{fibfibetale}) fail to be completely anisotropic, while the other 16 $\fib\boxtimes\fib$'s are completely anisotropic.

\subsubsection{$\mcal B\simeq psu(2)_7$}
The MFCs have four simple objects $\{1,X_{h_X},Y_{h_Y},Z_{h_Z}\}$ obeying monoidal products
\begin{table}[H]
\begin{center}
\begin{tabular}{c|c|c|c|c}
    $\otimes$&1&$X$&$Y$&$Z$\\\hline
    1&1&$X$&$Y$&$Z$\\\hline
    $X$&&$1\oplus Y$&$X\oplus Z$&$Y\oplus Z$\\\hline
    $Y$&&&$1\oplus Y\oplus Z$&$X\oplus Y\oplus Z$\\\hline
    $Z$&&&&$1\oplus X\oplus Y\oplus Z$
\end{tabular}.
\end{center}
\end{table}
\hspace{-17pt}Thus, they have
\[ \fp_{\mcal B}(1)=1,\quad\fp_{\mcal B}(X)=\frac{\sin\frac{2\pi}9}{\sin\frac\pi9},\quad\fp_{\mcal B}(Y)=\frac{\sin\frac{3\pi}9}{\sin\frac\pi9},\quad\fp_{\mcal B}(Z)=\frac{\sin\frac{4\pi}9}{\sin\frac\pi9}, \]
and
\[ \fp(\mcal B)=\frac9{4\sin^2\frac\pi9}\approx19.2. \]
The quantum dimensions are solutions of $d_X^2=1+d_Y,d_Xd_Y=d_X+d_Z,d_Xd_Z=d_Y+d_Z,d_Y^2=1+d_Y+d_Z,d_Yd_Z=d_X+d_Y+d_Z,d_Z^2=1+d_X+d_Y+d_Z$. There are three (non-zero) solutions
\[ (d_X,d_Y,d_Z)=(-\frac{\sin\frac\pi9}{\cos\frac\pi{18}},-\frac{\sin\frac\pi3}{\cos\frac\pi{18}},1-\frac{\sin\frac\pi9}{\cos\frac\pi{18}}),(-\frac{\cos\frac\pi{18}}{\sin\frac{2\pi}9},\frac{\sin\frac\pi3}{\sin\frac{2\pi}9},1-\frac{\cos\frac\pi{18}}{\sin\frac{2\pi}9}),(\frac{\sin\frac{2\pi}9}{\sin\frac\pi9},\frac{\sin\frac{3\pi}9}{\sin\frac\pi9},\frac{\sin\frac{4\pi}9}{\sin\frac\pi9}) \]
with categorical dimensions
\[ D^2(\mcal B)=\frac9{4\cos^2\frac\pi{18}}(\approx3.7),\quad\frac9{4\sin^2\frac{2\pi}9}(\approx5.4),\quad\frac9{4\sin^2\frac\pi9}, \]
respectively. They have conformal dimensions\footnote{To find these, study the necessary conditions $S_{i,j}\in\mbb R$ originating from self-duality of simple objects.}
\[ (h_X,h_Y,h_Z)=\begin{cases}(\frac13,\frac89,\frac23),(\frac23,\frac19,\frac13)&(\text{1st quantum dimensions}),\\(\frac13,\frac59,\frac23),(\frac23,\frac49,\frac13)&(\text{2nd quantum dimensions}),\\(\frac13,\frac29,\frac23),(\frac23,\frac79,\frac13)&(\text{3rd quantum dimensions}).\end{cases}\quad(\mods1) \]
Therefore, there are
\[ 3(\text{quantum dimensions})\times2(\text{conformal dimensions})\times2(\text{categorical dimensions})=12 \]
MFCs. Among them, those four with the third quantum dimensions give unitary MFCs, consistent with \cite{RSW07}. We classify connected étale algebras in all 12 MFCs simultaneously.

Just as in the previous example, we directly work with an ansatz
\[ A\cong1\oplus n_XX\oplus n_YY\oplus n_ZZ \]
with $n_X,n_Y,n_Z\in\mbb N$. It has
\[ \fp_{\mcal B}(A)=1+\frac1{\sin\frac\pi9}\left(n_X\sin\frac{2\pi}9+n_Y\sin\frac{3\pi}9+n_Z\sin\frac{4\pi}9\right). \]
For this to obey (\ref{fpdimAbound}), there are only four possibilities
\[ (n_X,n_Y,n_Z)=(0,0,0),(1,0,0),(0,1,0),(0,0,1). \]
The first candidate is nothing but the trivial connected étale algebra $A\cong1$ giving $\mcal B_A^0\simeq\mcal B_A\simeq\mcal B$. The last three candidates fail to be commutative because they do not satisfy the necessary condition (\ref{commutativenecessary}), and they are ruled out.

To summarize, we found
\begin{table}[H]
\begin{center}
\begin{tabular}{c|c|c|c}
    Connected étale algebra $A$&$\mcal B_A$&$\rank(\mcal B_A)$&Lagrangian?\\\hline
    1&$\mcal B$&4&No
\end{tabular}.
\end{center}
\caption{Connected étale algebras in rank four MFC $\mcal B\simeq psu(2)_7$}\label{rank4psu27results}
\end{table}
\hspace{-17pt}That is, all the 12 ambient MFCs $psu(2)_7$'s are completely anisotropic.

\subsection{Rank five}

\subsubsection{$\mcal B\simeq su(5)_1$}
The MFC has five simple objects $\{1,X_{h_X},Y_{h_Y},Z_{h_Z},W_{h_W}\}$ obeying monoidal products
\begin{table}[H]
\begin{center}
\begin{tabular}{c|c|c|c|c|c}
    $\otimes$&1&$X$&$Y$&$Z$&$W$\\\hline
    1&1&$X$&$Y$&$Z$&$W$\\\hline
    $X$&&$W$&$1$&$Y$&$Z$\\\hline
    $Y$&&&$Z$&$W$&$X$\\\hline
    $Z$&&&&$X$&$1$\\\hline
    $W$&&&&&$Y$
\end{tabular}.
\end{center}
\end{table}
\hspace{-17pt}Thus, they have
\[ \fp_{\mcal B}(1)=\fp_{\mcal B}(X)=\fp_{\mcal B}(Y)=\fp_{\mcal B}(Z)=\fp_{\mcal B}(W)=1, \]
and
\[ \fp(\mcal B)=5. \]
The quantum dimensions $d_X,d_Y,d_Z,d_W$ are solutions of $d_X^2=d_W,d_Xd_Y=1,d_Xd_Z=d_Y,d_Xd_W=d_Z,d_Y^2=d_Z,d_Yd_Z=d_W,d_Yd_W=d_X,d_Z^2=d_X,d_Zd_W=1,d_W^2=d_Y$. There is only one solution
\[ d_X=d_Y=d_Z=d_W=1, \]
and categorical dimension
\[ D^2(\mcal B)=5. \]
Thus, all MFCs are unitary.

Since the simple objects generate $\mbb Z/5\mbb Z$, we know their conformal dimensions take values in $\mbb Z/5$ mod 1. Concretely, we find two conformal dimensions\footnote{One of the easiest ways to find these values is to check $S^2=C$, where $C$ is the charge conjugation matrix. Naively, $(h_X,h_Y,h_Z,h_W)=(\frac35,\frac35,\frac25,\frac25),(\frac45,\frac45,\frac15,\frac15)$ also satisfy the conditions, but they are related to the two by changing names of simple objects.}
\[ (h_X,h_Y,h_Z,h_W)=(\frac15,\frac15,\frac45,\frac45),(\frac25,\frac25,\frac35,\frac35)\quad(\mods1). \]
Therefore, there are
\[ 1(\text{quantum dimension})\times2(\text{conformal dimensions})\times2(\text{categorical dimensions})=4 \]
MFCs. All of them are unitary.

The most general candidate for a connected algebra is given by
\[ A\cong1\oplus n_XX\oplus n_YY\oplus n_ZZ\oplus n_WW \]
with $n_j\in\mbb N$. It has
\[ \fp_{\mcal B}(A)=1+n_X+n_Y+n_Z+n_W. \]
For this to obey (\ref{fpdimAbound}), there are only five possibilities
\[ (n_X,n_Y,n_Z,n_W)=(0,0,0,0),(1,0,0,0),(0,1,0,0),(0,0,1,0),(0,0,0,1). \]
The first candidate is nothing but the trivial connected étale algebra $A\cong1$ giving $\mcal B_A^0\simeq\mcal B_A\simeq\mcal B$. The other four possibilities contain nontrivial simple object $b\not\cong1$. Since $b$ is not self-dual, the candidate $A\cong1\oplus b$ is not self-dual, leading to a contradiction. Thus, the four possibilities are ruled out. (One can also rule them out by realizing that they fail to be commutative.)

We conclude
\begin{table}[H]
\begin{center}
\begin{tabular}{c|c|c|c}
    Connected étale algebra $A$&$\mcal B_A$&$\rank(\mcal B_A)$&Lagrangian?\\\hline
    1&$\mcal B$&5&No
\end{tabular}.
\end{center}
\caption{Connected étale algebras in rank five MFC $\mcal B\simeq su(5)_1$}\label{rank5su51results}
\end{table}
\hspace{-17pt}That is, all four MFCs $su(5)_1$'s are completely anisotropic.

\subsubsection{$\mcal B\simeq su(2)_4$}
The MFC has five simple objects $\{1,X_{h_X},Y_{h_Y},Z_{h_Z},W_{h_W}\}$ obeying monoidal products
\begin{table}[H]
\begin{center}
\begin{tabular}{c|c|c|c|c|c}
    $\otimes$&1&$X$&$Y$&$Z$&$W$\\\hline
    1&1&$X$&$Y$&$Z$&$W$\\\hline
    $X$&&$1$&$Z$&$Y$&$W$\\\hline
    $Y$&&&$1\oplus W$&$X\oplus W$&$Y\oplus Z$\\\hline
    $Z$&&&&$1\oplus W$&$Y\oplus Z$\\\hline
    $W$&&&&&$1\oplus X\oplus W$
\end{tabular}.
\end{center}
\end{table}
\hspace{-17pt}Thus, they have
\[ \fp_{\mcal B}(1)=1=\fp_{\mcal B}(X),\quad\fp_{\mcal B}(Y)=\sqrt3=\fp_{\mcal B}(Z),\quad\fp_{\mcal B}(W)=2, \]
and
\[ \fp(\mcal B)=12. \]
The quantum dimensions are solutions of $d_X^2=1,d_Xd_Y=d_Z,d_Xd_Z=d_Y,d_Xd_W=d_W,d_Y^2=1+d_W,d_Yd_Z=d_X+d_W,d_Yd_W=d_Y+d_Z,d_Z^2=1+d_W,d_Zd_W=d_Y+d_Z,d_W^2=1+d_X+d_W$. There are two solutions
\[ (d_X,d_Y,d_Z,d_W)=(1,-\sqrt3,-\sqrt3,2),(1,\sqrt3,\sqrt3,2). \]
They both have the same categorical dimension
\[ D^2(\mcal B)=12. \]

The conformal dimensions can be computed studying necessary conditions such as $S_{i,j}\in\mbb R$ or non-degeneracy of $S$-matrix. As a result, independent of quantum dimensions, one finds four conformal dimensions\footnote{Naively, one can also have
\[ (h_X,h_Y,h_Z,h_W)=(0,\frac58,\frac18,\pm\frac13),(0,\frac78,\frac38,\pm\frac13)\quad(\mods1), \]
but these give the same MFCs as those in the main text under change of names $Y\leftrightarrow Z$.}
\[ (h_X,h_Y,h_Z,h_W)=(0,\frac18,\frac58,\frac13),(0,\frac18,\frac58,\frac23),(0,\frac38,\frac78,\frac13),(0,\frac38,\frac78,\frac23)\quad(\mods1). \]
Therefore, there are
\[ 2(\text{quantum dimensions})\times4(\text{conformal dimensions})\times2(\text{categorical dimensions})=16 \]
MFCs, consistent with \cite{BNRW15}. Among them, those eight with positive quantum dimensions give unitary MFCs. We study connected étale algebras in all 16 MFCs simultaneously.

The most general candidate for a connected algebra is given by
\[ A\cong1\oplus n_XX\oplus n_YY\oplus n_ZZ\oplus n_WW \]
with $n_j\in\mbb N$. It has
\[ \fp_{\mcal B}(A)=1+n_X+\sqrt3(n_Y+n_Z)+2n_W. \]
For this to obey (\ref{fpdimAbound}), there are only six possibilities
\[ (n_X,n_Y,n_Z,n_W)=(0,0,0,0),(1,0,0,0),(2,0,0,0),(0,1,0,0),(0,0,1,0),(0,0,0,1). \]
The first candidate is nothing but the trivial connected étale algebra $A\cong1$ giving $\mcal B_A^0\simeq\mcal B_A\simeq\mcal B$. The second corresponds to the $\mbb Z/2\mbb Z$ algebra $A\cong1\oplus X$. Since $X$ has $(d_X,h_X)=(1,0)$, we know \cite{KK23etale} $c_{X,X}\cong id_1$, and the candidate does give connected étale algebra. Let us study what is the category of right $A$-modules. The connected étale algebra has
\[ \fp_{\mcal B}(A)=2, \]
and (\ref{FPdimBA0}) demands
\[ \fp(\mcal B_A^0)=3 \]
with
\[ c(\mcal B_A^0)=c(\mcal B)=2\quad(\mods4). \]
The only possibility is
\[ \mcal B_A^0\simeq\vecG_{\mbb Z/3\mbb Z}^1. \]
The category of right $A$-modules $\mcal B_A$ contains this MFC as a subcategory. It has
\[ \fp(\mcal B_A)=\frac{\fp(\mcal B)}{\fp_{\mcal B}(A)}=6. \]
There are a few candidates. For example, multiplicity-free fusion categories $\text{TY}(\mbb Z/3\mbb Z),D_3,\vecG_{\mbb Z/6\mbb Z}^\alpha$ with $\fp=6$ contain $\mcal B_A^0\simeq\vecG_{\mbb Z/3\mbb Z}^1$, and they can be $\mcal B_A$. It turns out that the correct one is
\begin{equation}
    \mcal B_A\simeq\text{TY}(\mbb Z/3\mbb Z),\label{rank5su24BA}
\end{equation}
the $\mbb Z/3\mbb Z$ Tambara-Yamagami category \cite{TY98}. When $\mcal B\simeq\mcal C(A_1,4)$, this result is well-known; the connected étale algebra gives $\mcal B_A$ corresponding to the $D_4$ Dynkin diagram \cite{KO01}, and especially the category has $\rank(\mcal B_A)=4$. For the other MFCs, one of the easiest ways to find this fact is to perform anyon condensation \cite{BSS02,BSS02'}. It is a process to get a new MFC $\mcal B_A^0$ from an old MFC $\mcal B$ for an algebra $A\in\mcal B$. More concretely, one `identifies' condensing object ($X$ in our case) with the unit $1\in\mcal B$. Since $Z\cong X\otimes Y,Y\cong X\otimes Z$, this `identifies' $Y$ and $Z$. Furthermore, $W$ with quantum dimension two `splits' into two invertible simple objects \cite{BS08}. As a result, one obtains\footnote{More mathematically, we should search for NIM-reps. We carried out the computation, and found a four-dimensional NIM rep
\[ n_X=1_4,\quad n_Y=\begin{pmatrix}0&0&0&1\\0&0&0&1\\0&0&0&1\\1&1&1&0\end{pmatrix}=n_Z,\quad n_W=\begin{pmatrix}0&1&1&0\\1&0&1&0\\1&1&0&0\\0&0&0&2\end{pmatrix}. \]
Denoting a basis of $\mcal M$ by $\{m_1,m_2,m_3,m_4\}$, we get a multiplication table
\begin{table}[H]
\begin{center}
\begin{tabular}{c|c|c|c|c}
    $b_j\otimes\backslash$&$F(m_1)$&$F(m_2)$&$F(m_3)$&$F(m_4)$\\\hline
    1&$F(m_1)$&$F(m_2)$&$F(m_3)$&$F(m_4)$\\
    $X$&$F(m_1)$&$F(m_2)$&$F(m_3)$&$F(m_4)$\\
    $Y$&$F(m_4)$&$F(m_4)$&$F(m_4)$&$F(m_1)\oplus F(m_2)\oplus F(m_3)$\\
    $Z$&$F(m_4)$&$F(m_4)$&$F(m_4)$&$F(m_1)\oplus F(m_2)\oplus F(m_3)$\\
    $W$&$F(m_2)\oplus F(m_3)$&$F(m_1)\oplus F(m_3)$&$F(m_1)\oplus F(m_2)$&$2F(m_4)$
\end{tabular}.
\end{center}
\end{table}
\hspace{-14pt}
Here, for a fixed $m$ in a left $\mcal B$-module category $\mcal M$, $F:=\Hom(m,-):\mcal M\to\mcal B_{\Hom(m,m)}$ is a functor to the category of right $\Hom(m,m)$-modules. The $\Hom(m_1,m_2)\in\mcal C$ is the internal Hom \cite{O01} defined by natural isomorphism
\[ m_1,m_2\in\mcal M,\forall c\in\mcal C,\quad\mcal M(c\triangleright m_1,m_2)\cong\mcal C(c,\Hom(m_1,m_2)). \]
For every $m\in\mcal M$, $\Hom(m,m)$ has a canonical structure of an algebra in $\mcal C$.

In a suitable basis, we can identify
\[ F(m_1)\cong1\oplus X,\quad F(m_2)\cong W\cong F(m_3),\quad F(m_4)\cong Y\oplus Z. \]
They have (\ref{dBAm})
\[ d_{\mcal B_A}(F(m_1))=d_{\mcal B_A}(F(m_2))=d_{\mcal B_A}(F(m_3))=1,\quad d_{\mcal B_A}(F(m_4))=\pm\sqrt3, \]
showing $\mcal B_A\simeq\text{TY}(\mbb Z/3\mbb Z)$, especially $\rank(\mcal B_A)=4$.} a $\text{TY}(\mbb Z/3\mbb Z)$. (The $\mcal B_A$ is called a broken phase.)

The third candidate $A\cong1\oplus2X$ has $\fp_{\mcal B}(A)=3$, and demands $\fp(\mcal B_A^0)=\frac43$. Since there is no MFC with this Frobenius-Perron dimension, the candidate is ruled out. The last three possibilities contain nontrivial simple object $b\not\cong1$ with nontrivial conformal dimension. Thus, they fail to be commutative, and are ruled out.

To sum up, we find
\begin{table}[H]
\begin{center}
\begin{tabular}{c|c|c|c}
    Connected étale algebra $A$&$\mcal B_A$&$\rank(\mcal B_A)$&Lagrangian?\\\hline
    1&$\mcal B$&5&No\\
    $1\oplus X$&$\text{TY}(\mbb Z/3\mbb Z)$&4&No
\end{tabular}.
\end{center}
\caption{Connected étale algebras in rank five MFC $\mcal B\simeq su(2)_4$}\label{rank5su24results}
\end{table}
\hspace{-17pt}Therefore, all 16 $\mcal B\simeq su(2)_4$'s fail to be completely anisotropic.

\subsubsection{$\mcal B\simeq psu(2)_9$}
The MFC has five simple objects $\{1,X_{h_X},Y_{h_Y},Z_{h_Z},W_{h_W}\}$ obeying monoidal products
\begin{table}[H]
\begin{center}
\begin{tabular}{c|c|c|c|c|c}
    $\otimes$&1&$X$&$Y$&$Z$&$W$\\\hline
    1&1&$X$&$Y$&$Z$&$W$\\\hline
    $X$&&$1\oplus Y$&$X\oplus Z$&$Y\oplus W$&$Z\oplus W$\\\hline
    $Y$&&&$1\oplus Y\oplus W$&$X\oplus Z\oplus W$&$Y\oplus Z\oplus W$\\\hline
    $Z$&&&&$1\oplus Y\oplus Z\oplus W$&$X\oplus Y\oplus Z\oplus W$\\\hline
    $W$&&&&&$1\oplus X\oplus Y\oplus Z\oplus W$
\end{tabular}.
\end{center}
\end{table}
\hspace{-17pt}Thus, they have
\[ \hspace{-50pt}\fp_{\mcal B}(1)=1,\quad\fp_{\mcal B}(X)=\frac{\sin\frac{2\pi}{11}}{\sin\frac\pi{11}},\quad\fp_{\mcal B}(Y)=\frac{\sin\frac{3\pi}{11}}{\sin\frac\pi{11}},\quad\fp_{\mcal B}(Z)=\frac{\sin\frac{4\pi}{11}}{\sin\frac\pi{11}},\quad\fp_{\mcal B}(W)=\frac{\sin\frac{5\pi}{11}}{\sin\frac\pi{11}}, \]
and
\[ \fp(\mcal B)=\frac{11}{4\sin^2\frac\pi{11}}\approx34.6. \]
The quantum dimensions are solutions of $d_X^2=1+d_Y,d_Xd_Y=d_X+d_Z,d_Xd_Z=d_Y+d_W,d_Xd_W=d_Z+d_W,d_Y^2=1+d_Y+d_W,d_Yd_Z=d_X+d_Z+d_W,d_Yd_W=d_Y+d_Z+d_W,d_Z^2=1+d_Y+d_Z+d_W,d_Zd_W=d_X+d_Y+d_Z+d_W,d_W^2=1+d_X+d_Y+d_Z+d_W$. There are five solutions
\begin{align*}
    (d_X,d_Y,d_Z,d_W)=&(\frac{\sin\frac\pi{11}}{\cos\frac\pi{22}},-\frac{\sin\frac{4\pi}{11}}{\cos\frac\pi{22}},-\frac{\sin\frac{2\pi}{11}}{\cos\frac\pi{22}},\frac{\sin\frac{3\pi}{11}}{\cos\frac\pi{22}}),(-\frac{\sin\frac{3\pi}{11}}{\cos\frac{3\pi}{22}},-\frac{\sin\frac\pi{11}}{\cos\frac{3\pi}{22}},\frac{\sin\frac{5\pi}{11}}{\cos\frac{3\pi}{22}},-\frac{\sin\frac{2\pi}{11}}{\cos\frac{3\pi}{22}}),\\
    &(\frac{\sin\frac{5\pi}{11}}{\sin\frac{3\pi}{11}},\frac{\sin\frac{2\pi}{11}}{\sin\frac{3\pi}{11}},-\frac{\sin\frac\pi{11}}{\sin\frac{3\pi}{11}},-\frac{\sin\frac{4\pi}{11}}{\sin\frac{3\pi}{11}}),(-\frac{\sin\frac{4\pi}{11}}{\sin\frac{2\pi}{11}},\frac{\sin\frac{5\pi}{11}}{\sin\frac{2\pi}{11}},-\frac{\sin\frac{3\pi}{11}}{\sin\frac{2\pi}{11}},\frac{\sin\frac\pi{11}}{\sin\frac{2\pi}{11}}),\\
    &(\frac{\sin\frac{2\pi}{11}}{\sin\frac\pi{11}},\frac{\sin\frac{3\pi}{11}}{\sin\frac\pi{11}},\frac{\sin\frac{4\pi}{11}}{\sin\frac\pi{11}},\frac{\sin\frac{5\pi}{11}}{\sin\frac\pi{11}}).
\end{align*}
Only the last quantum dimensions give unitary MFCs. They have categorical dimensions
\[ D^2(\mcal B)=\frac{11}{4\cos^2\frac\pi{22}}(\approx2.8),\quad\frac{11}{4\cos^2\frac{3\pi}{22}}(\approx3.3),\quad\frac{11}{4\sin^2\frac{3\pi}{11}}(\approx4.8),\quad\frac{11}{4\sin^2\frac{2\pi}{11}}(\approx9.4),\quad\frac{11}{4\sin^2\frac\pi{11}}, \]
respectively. Studying necessary conditions, one finds their conformal dimensions
\[ (h_X,h_Y,h_Z,h_W)=\begin{cases}
(\frac1{11},\frac{10}{11},\frac5{11},\frac8{11}),(\frac{10}{11},\frac1{11},\frac6{11},\frac3{11})&(\text{1st quantum dimensions}),\\
(\frac3{11},\frac8{11},\frac4{11},\frac2{11}),(\frac8{11},\frac3{11},\frac7{11},\frac9{11})&(\text{2nd quantum dimensions}),\\
(\frac5{11},\frac6{11},\frac3{11},\frac7{11}),(\frac6{11},\frac5{11},\frac8{11},\frac4{11})&(\text{3rd quantum dimensions}),\\
(\frac4{11},\frac7{11},\frac9{11},\frac{10}{11}),(\frac7{11},\frac4{11},\frac2{11},\frac1{11})&(\text{4th quantum dimensions}),\\
(\frac2{11},\frac9{11},\frac{10}{11},\frac5{11}),(\frac9{11},\frac2{11},\frac1{11},\frac6{11}),&(\text{5th quantum dimensions}).\end{cases}\quad(\mods1) \]
Therefore, we find
\[ 5(\text{quantum dimensions})\times2(\text{conformal dimensions})\times2(\text{categorical dimensions})=20 \]
MFCs, among which those four with the fifth quantum dimensions give unitary MFCs. We study connected étale algebras in all 20 MFCs simultaneously.

The most general candidate for a connected algebra is given by
\[ A\cong1\oplus n_XX\oplus n_YY\oplus n_ZZ\oplus n_WW \]
with $n_j\in\mbb N$. It has
\[ \fp_{\mcal B}(A)=1+\frac1{\sin\frac\pi{11}}\left(n_X\sin\frac{2\pi}{11}+n_Y\sin\frac{3\pi}{11}+n_Z\sin\frac{4\pi}{11}+n_W\sin\frac{5\pi}{11}\right). \]
For this to obey (\ref{fpdimAbound}), there are only seven possibilities
\[ (n_X,n_Y,n_Z,n_W)=(0,0,0,0),(1,0,0,0),(2,0,0,0),(0,1,0,0),(1,1,0,0),(0,0,1,0),(0,0,0,1). \]
The first candidate is nothing but the trivial connected étale algebra $A\cong1$ giving $\mcal B_A^0\simeq\mcal B_A\simeq\mcal B$. The other six possibilities contain nontrivial simple object(s) $b\not\cong1$. Thus, they do not give connected étale algebras because they fail to satisfy the necessary condition.

To summarize, we find
\begin{table}[H]
\begin{center}
\begin{tabular}{c|c|c|c}
    Connected étale algebra $A$&$\mcal B_A$&$\rank(\mcal B_A)$&Lagrangian?\\\hline
    1&$\mcal B$&5&No
\end{tabular}.
\end{center}
\caption{Connected étale algebras in rank five MFC $\mcal B\simeq psu(2)_9$}\label{rank5psu29results}
\end{table}
\hspace{-17pt}Namely, all 20 $\mcal B\simeq psu(2)_9$'s are completely anisotropic.

\subsubsection{$\mcal B\simeq psu(3)_4$}
The MFC has five simple objects $\{1,X_{h_X},Y_{h_Y},Z_{h_Z},W_{h_W}\}$ obeying monoidal products\footnote{In terms of Dynkin labels in the $\wf{su}(3)_4$ WZW model, we have correspondences
\[ (0,0)\leftrightarrow1,\quad(1,1)\leftrightarrow Y,\quad(3,0)\leftrightarrow Z\ (\text{or }W),\quad(0,3)\leftrightarrow W\ (\text{or }Z),\quad(2,2)\leftrightarrow X. \]
In terms of \cite{GK94}, we have correspondences
\[ \phi_1\leftrightarrow X,\quad\phi_2\leftrightarrow Y,\quad\phi_3\leftrightarrow Z(\text{or }W),\quad\phi_4\leftrightarrow W(\text{or }Z). \]
Since fusion matrices are non-symmetric, we wrote all monoidal products. It turns out that they are symmetric. The $S$-matrices
\[ \widetilde S=\begin{pmatrix}1&d_X&d_Y&d_Z&d_W\\d_X&d_Y&-1&-d_Z&-d_W\\d_Y&-1&-d_X&d_Z&d_W\\d_Z&-d_Z&d_Z&\frac{-1\pm\sqrt7i}2d_Z&\frac{-1\mp\sqrt7i}2d_Z\\d_W&-d_W&d_W&\frac{-1\mp\sqrt7i}2d_Z&\frac{-1\pm\sqrt7i}2d_Z\end{pmatrix} \]
are also symmetric.}
\begin{table}[H]
\begin{center}
\begin{tabular}{c|c|c|c|c|c}
    $\otimes$&1&$X$&$Y$&$Z$&$W$\\\hline
    1&1&$X$&$Y$&$Z$&$W$\\\hline
    $X$&$X$&$1\oplus X\oplus Y$&$X\oplus Y\oplus Z\oplus W$&$Y\oplus W$&$Y\oplus Z$\\\hline
    $Y$&$Y$&$X\oplus Y\oplus Z\oplus W$&$1\oplus X\oplus2Y\oplus Z\oplus W$&$X\oplus Y\oplus Z$&$X\oplus Y\oplus W$\\\hline
    $Z$&$Z$&$Y\oplus W$&$X\oplus Y\oplus Z$&$X\oplus W$&$1\oplus Y$\\\hline
    $W$&$W$&$Y\oplus Z$&$X\oplus Y\oplus W$&$1\oplus Y$&$X\oplus Z$
\end{tabular}.
\end{center}
\end{table}
\hspace{-17pt}Note that this is the first example with multiplicity ${N_{Y,Y}}^Y=2$. Accordingly, fusion matrices are not symmetric:
\[ \hspace{-40pt}N_X=\begin{pmatrix}0&1&0&0&0\\1&1&1&0&0\\0&1&1&1&1\\0&0&1&0&1\\0&0&1&1&0\end{pmatrix},\quad N_Y=\begin{pmatrix}0&0&1&0&0\\0&1&1&1&1\\1&1&2&1&1\\0&1&1&1&0\\0&1&1&0&1\end{pmatrix},\quad N_Z=\begin{pmatrix}0&0&0&1&0\\0&0&1&0&1\\0&1&1&1&0\\0&1&0&0&1\\1&0&1&0&0\end{pmatrix},\quad N_W=\begin{pmatrix}0&0&0&0&1\\0&0&1&1&0\\0&1&1&0&1\\1&0&1&0&0\\0&1&0&1&0\end{pmatrix}. \]
They have
\[ \hspace{-30pt}\fp_{\mcal B}(1)=1,\quad\fp_{\mcal B}(X)=\frac{\sin\frac{3\pi}{14}}{\sin\frac\pi{14}},\quad\fp_{\mcal B}(Y)=\frac{\sin\frac{5\pi}{14}}{\sin\frac\pi{14}},\quad\fp_{\mcal B}(Z)=\frac1{2\sin\frac\pi{14}}=\fp_{\mcal B}(W), \]
and
\[ \fp(\mcal B)=\frac7{4\sin^2\frac\pi{14}}\approx35.3. \]
The quantum dimensions are solutions of $d_X^2=1+d_X+d_Y,d_Xd_Y=d_X+d_Y+d_Z+d_W,d_Xd_Z=d_Y+d_W,d_Xd_W=d_Y+d_Z,d_Y^2=1+d_X+2d_Y+d_Z+d_W,d_Yd_Z=d_X+d_Y+d_Z,d_Yd_W=d_X+d_Y+d_W,d_Z^2=d_X+d_W,d_Zd_W=1+d_Y,d_W^2=d_X+d_Z$. There are three solutions
\begin{align*}
    (d_X,d_Y,d_Z,d_W)=&(-\frac{\sin\frac\pi{14}}{\sin\frac{5\pi}{14}},-\frac{\sin\frac{3\pi}{14}}{\sin\frac{5\pi}{14}},\frac1{2\sin\frac{5\pi}{14}},\frac1{2\sin\frac{5\pi}{14}}),(\frac{\sin\frac{5\pi}{14}}{\cos\frac{2\pi}7},-\frac{\sin\frac\pi{14}}{\cos\frac{2\pi}7},-\frac1{2\cos\frac{2\pi}7},-\frac1{2\cos\frac{2\pi}7}),\\
    &(\frac{\sin\frac{3\pi}{14}}{\sin\frac\pi{14}},\frac{\sin\frac{5\pi}{14}}{\sin\frac\pi{14}},\frac1{2\sin\frac\pi{14}},\frac1{2\sin\frac\pi{14}}).
\end{align*}
Only the last quantum dimensions give unitary MFCs. They have categorical dimensions
\[ D^2(\mcal B)=\frac7{4\sin^2\frac{5\pi}{14}}(\approx2.2),\quad\frac7{4\cos^2\frac{2\pi}7}(\approx4.5),\quad\frac7{4\sin^2\frac\pi{14}}, \]
respectively. Studying necessary conditions, one finds they have conformal dimensions
\[ (h_X,h_Y,h_Z,h_W)=\begin{cases}(\frac27,\frac67,\frac57,\frac57),(\frac57,\frac17,\frac27,\frac27)&(\text{1st quantum dimensions}),\\(\frac37,\frac27,\frac47,\frac47),(\frac47,\frac57,\frac37,\frac37)&(\text{2nd quantum dimensions}),\\(\frac17,\frac37,\frac67,\frac67),(\frac67,\frac47,\frac17,\frac17)&(\text{3rd quantum dimensions}).\end{cases}\quad(\mods1) \]
Therefore, there are
\[ 3(\text{quantum dimensions})\times2(\text{conformal dimensions})\times2(\text{categorical dimensions})=12 \]
MFCs, among which those four with the last quantum dimensions give unitary MFCs.

Having found the conformal dimensions, we can classify connected étale algebras. The most general ansatz for a connected algebra is given by
\[ A\cong1\oplus n_XX\oplus n_YY\oplus n_ZZ\oplus n_WW \]
with $n_j\in\mbb N$. It has
\[ \fp_{\mcal B}(A)=1+\frac1{\sin\frac\pi{14}}\left(n_X\sin\frac{3\pi}{14}+n_Y\sin\frac{5\pi}{14}+\frac12n_Z+\frac12n_W\right). \]
For this to obey (\ref{fpdimAbound}), there are only eight possibilities
\[ (n_X,n_Y,n_Z,n_W)=(0,0,0,0),(1,0,0,0),(0,1,0,0),(0,0,1,0),(0,0,2,0),(0,0,1,1),(0,0,0,1),(0,0,0,2). \]
The first candidate is nothing but the trivial connected étale algebra $A\cong1$ giving $\mcal B_A^0\simeq\mcal B_A\simeq\mcal B$. The other seven possibilities contain nontrivial simple object(s) $b\not\cong1$. Thus, they fail to be commutative because they do not satisfy the necessary condition (\ref{commutativenecessary}).

To summarize, we find
\begin{table}[H]
\begin{center}
\begin{tabular}{c|c|c|c}
    Connected étale algebra $A$&$\mcal B_A$&$\rank(\mcal B_A)$&Lagrangian?\\\hline
    1&$\mcal B$&5&No
\end{tabular}.
\end{center}
\caption{Connected étale algebras in rank five MFC $\mcal B\simeq psu(3)_4$}\label{rank5psu34results}
\end{table}
\hspace{-17pt}Namely, all 12 MFCs $\mcal B\simeq psu(3)_4$'s are completely anisotropic.

\section{Physical applications}
Our classification results have physical implications on ground state degeneracy (GSD) and spontaneous symmetry breaking (SSB). In this section, we discuss the physical applications.

\subsection{Theorems}
Let $\mcal C$ be a fusion category. Two-dimensional gapped phases with $\mcal C$ symmetry stand in bijection with $\mcal C$-module categories \cite{TW19,HLS21}
\begin{equation}
    \{\text{2d }\mcal C\text{-symmetric gapped phases}\}\cong\{\mcal C\text{-module categories }\mcal M\}.\label{1to1}
\end{equation}
In particular, GSD in the LHS is equal to the $\rank(\mcal M)$ of a $\mcal C$-module category
\[ \text{GSD}=\rank(\mcal M). \]
This leads to the\newline

\textbf{Theorem.} \textit{Let $\mcal B$ be a modular fusion category up to rank five and $A\in\mcal B$ a connected étale algebra. Suppose two-dimensional $\mcal B$-symmetric gapped phases are described by indecomposable $\mcal B_A$'s. Then, the gapped phases have}
\[ \text{GSD}\in\begin{cases}\{1\}&(\mcal B\simeq\vect),\\\{2\}&(\mcal B\simeq\vecG_{\mbb Z/2\mbb Z}^{-1}),\\\{2\}&(\mcal B\simeq\fib),\\\{3\}&(\mcal B\simeq\vecG_{\mbb Z/3\mbb Z}^1),\\\{3\}&(\mcal B\simeq\ising),\\\{3\}&(\mcal B\simeq psu(2)_5),\\\{2,4\}&(\mcal B\simeq\vecG_{\mbb Z/2\mbb Z}^{-1}\boxtimes\vecG_{\mbb Z/2\mbb Z}^{-1}\ (\text{four with }(\ref{Z2Z2etale}))),\\\{4\}&(\mcal B\simeq\vecG_{\mbb Z/2\mbb Z}^{-1}\boxtimes\vecG_{\mbb Z/2\mbb Z}^{-1}\ (\text{the other 16})),\\\{2,4\}&(\mcal B\simeq\tc\ (\text{four with }(\ref{rank4ToricCoderesultalg}))),\\\{4\}&(\mcal B\simeq\tc\text{ (the other six)}),\\\{4\}&(\mcal B\simeq\vecG_{\mbb Z/4\mbb Z}^\alpha),\\\{4\}&(\mcal B\simeq su(2)_3),\\\{2,4\}&(\mcal B\simeq\fib\boxtimes\fib\ (\text{four with }(\ref{fibfibetale}))),\\\{4\}&(\mcal B\simeq\fib\boxtimes\fib\ (\text{the other 16})),\\\{4\}&(\mcal B\simeq psu(2)_7),\\\{5\}&(\mcal B\simeq\vecG_{\mbb Z/5\mbb Z}^1),\\\{4,5\}&(\mcal B\simeq su(2)_4),\\\{5\}&(\mcal B\simeq psu(2)_9),\\\{5\}&(\mcal B\simeq psu(3)_4).\end{cases} \]\newline

The Theorem also proves certain SSBs. Here, we have the\newline

\textbf{Definition.} \cite{KK23GSD} Let $\mcal C$ be a fusion category and $\mcal M$ be a (left) $\mcal C$-module category describing a $\mcal C$-symmetric gapped phase. A symmetry $c\in\mcal C$ is called \textit{spontaneously broken} if $\exists m\in\mcal M$ such that $c\triangleright m\not\cong m$. We also say $\mcal C$\textit{ is spontaneously broken} if there exists a spontaneously broken object $c\in\mcal C$. A categorical symmetry $\mcal C$ is called \textit{preserved} (i.e., not spontaneously broken) if all objects act trivially.\newline

With the definition, one can show a\newline

\textbf{Lemma.} \cite{KK23GSD} \textit{Let $\mcal C$ be a fusion category and $\mcal M$ be an indecomposable (left) $\mcal C$-module category. Then, $\rank(\mcal M)>1$ implies SSB of $\mcal C$ (i.e., $\mcal C$ is spontaneously broken).}\newline

Therefore, we have proved SSBs:\newline

\textbf{Theorem.} \textit{Let $\mcal B$ be a modular fusion category and $A\in\mcal B$ be a connected étale algebra. In two-dimensional $\mcal B$-symmetric gapped phases described by indecomposable $\mcal B_A$'s, $\mcal B$ symmetries are spontaneously broken for}
\[ \mcal B\simeq\begin{cases}\vecG_{\mbb Z/2\mbb Z}^{-1},\\\fib,\\\vecG_{\mbb Z/3\mbb Z}^1,\\\ising,\\psu(2)_5,\\\vecG_{\mbb Z/2\mbb Z}^{-1}\boxtimes\vecG_{\mbb Z/2\mbb Z}^{-1},\\\tc,\\\vecG_{\mbb Z/4\mbb Z}^\alpha,\\su(2)_3,\\\fib\boxtimes\fib,\\psu(2)_7,\\\vecG_{\mbb Z/5\mbb Z}^1,\\su(2)_4,\\psu(2)_9,\\psu(3)_4.\end{cases} \]
\textit{Namely, all modular fusion categories up to rank five but $\mcal B\simeq\vect$ are spontaneously broken.}\newline

\textbf{Remark.} As noted in \cite{KK23GSD}, commutativity of an algebra seems too strong; numerical computation suggests an existence of $\mcal B$-symmetric gapped phase described by $\mcal B_A$ with non-commutative connected separable algebra.\newline

This result motivates the\newline

\textbf{Theorem.} \textit{Let $\mc B$ be a modular fusoin category with $\rank(\mc B)>1$. The $\mc B$ symmetry is spontaneously broken in two-dimensional $\mc B$-symmetric gapped phases described by indecomposable $\mc B_A$'s with connected étale algebra $A\in\mc B$.}\newline

\textit{Proof.} Assume the opposite. The contrapositive of the lemma claims $\rank(\mc B_A)=1$, or $\mc B_A\simeq\vect$. The formula $\fp(\mc B_A)=\frac{\fp(\mc B)}{\fp_{\mc B}(A)}$ together with $\fp(\vect)=1$ demands $\fp_{\mc B}(A)=\fp(\mc B)$, but there does not exist an MFC with $\fp(\mc B_A^0)=\frac{\fp(\mc B)}{(\fp_{\mc B}(A))^2}=\frac1{\fp(\mc B)}$ for $\rank(\mc B)>1$. $\square$\newline

This explains, assuming the gapped phases are described by such module categories, why massive deformations of unitary minimal models\footnote{We basically follow the notations of \cite{FMS}.} have degenerated vacua. Concretely, a relevant operator $\phi_{1,3}$ in the unitary minimal model $M(2M+3,2M+2)$ with $M=1,2,3,\dots$ preserves rank $(2M+1)$ MFC $\mc B$, and negative Lagrangian coupling leads to a gapped phase with $\text{GSD}=(2M+1)$ \cite{KK21}. The GSD is understood as a consequence of spontaneous $\mc B$-symmetry breaking.

\subsection{Examples}
Let us discuss some concrete examples and \textit{predict} GSDs. Since unitary cases are relatively well-understood, we study less-understood non-unitary examples.

Pick a non-unitary minimal model $M(p,2p+1)$ with $p\ge2$. It was proved \cite{KK22free} that its relevant $\phi_{5,1}$-deformation preserves modular fusion subcategory. If the relevant deformation triggers massless renormalization group (RG) flow, it is known \cite{Z90,Z91,M91,RST94,DDT00} that the IR theory is another non-unitary minimal model $M(p,2p-1)$. Its further relevant $\phi_{1,2}$-deformation also preserves modular fusion subcategory. Typically, massless RG flows require fine-tuning. For example, in infinitely many massless RG flows among minimal models $M(m+1,m)\to M(m,m-1)$, one has to choose the correct signs of relevant couplings. The other signs generally trigger massive RG flows to gapped phases. Therefore, relevant deformations of either $M(p,2p+1)$ or $M(p,2p-1)$ would generally lead to $\mcal B$-symmetric gapped phases with modular $\mcal B$. In this situation, we can apply our classification results. Below, we study RG flows triggered by relevant deformations of the non-unitary minimal models.\newline

\paragraph{$M(5,9)+\phi_{1,2}$.} The relevant deformation preserves rank four MFC $\mcal B$ with simple objects \cite{KK22free} $\{\mcal L_{1,1},\mcal L_{3,1},\mcal L_{5,1},\mcal L_{7,1}\}$. They form $\mcal B\simeq psu(2)_7$ with identifications
\[ X\cong\mcal L_{7,1},\quad Y\cong\mcal L_{3,1},\quad Z\cong\mcal L_{5,1}. \]
They have the first (non-unitary) quantum dimensions. Their conformal dimensions
\[ (h_{7,1},h_{3,1},h_{5,1})=(\frac{11}3,\frac19,\frac43) \]
match our second conformal dimensions mod 1. Thus, our classification result immediately implies that the massive RG flow described by $\mcal B_A$ has $\text{GSD}=4$ and $\mcal B$ symmetry is spontaneously broken.

\paragraph{$M(5,11)+\phi_{5,1}$.} The relevant deformation preserves rank four MFC $\mcal B$ with simple objects \cite{KK22free} $\{\mcal L_{1,1},\mcal L_{1,2},\mcal L_{1,3},\mcal L_{1,4}\}$. They form $\mcal B\simeq su(2)_3$ with identifications
\[ X\cong\mcal L_{1,4},\quad Y\cong\mcal L_{1,3},\quad Z\cong\mcal L_{1,2}. \]
They have the first quantum dimensions. (Thus, the MFC is actually unitary.) Their conformal dimensions
\[ (h_{1,4},h_{1,3},h_{1,2})=(\frac{27}4,\frac{17}5,\frac{23}{20}) \]
match our third conformal dimensions mod 1. Having specified the symmetry, we immediately learn the massive RG flow described by $\mcal B_A$ has $\text{GSD}=4$ and $\mcal B$ symmetry is spontaneously broken.

\paragraph{$M(6,11)+\phi_{1,2}$.} The relevant deformation preserves rank five MFC $\mcal B$ with simple objects \cite{KK22free} $\{\mcal L_{1,1},\mcal L_{3,1},\mcal L_{5,1},\mcal L_{7,1},\mcal L_{9,1}\}$. They form $\mcal B\simeq psu(2)_9$ with identifications
\[ X\cong\mcal L_{9,1},\quad Y\cong\mcal L_{3,1},\quad Z\cong\mcal L_{7,1},\quad W\cong\mcal L_{5,1}. \]
They have the first (non-unitary) quantum dimensions. Their conformal dimensions
\[ (h_{9,1},h_{3,1},h_{7,1},h_{5,1})=(\frac{76}{11},\frac1{11},\frac{39}{11},\frac{14}{11}) \]
match our second conformal dimensions mod 1. Thus, our classification result immediately implies the following. The massive RG flow described by $\mcal B_A$ has $\text{GSD}=5$ and $\mcal B$ symmetry is spontaneously broken.

\paragraph{$M(6,13)+\phi_{5,1}$.} The relevant deformation preserves rank five MFC with simple objects \cite{KK22free} $\{\mcal L_{1,1},\mcal L_{1,2},\mcal L_{1,3},\mcal L_{1,4},\mcal L_{1,5}\}$. They form $\mcal B\simeq su(2)_4$ with identifications\footnote{Just from monoidal products or quantum dimensions, we cannot fix the $\mbb Z/2\mbb Z$ ambiguity $Y\leftrightarrow Z$, but conformal dimensions fix the ambiguity and lead to the unique identifications.}
\[ X\cong\mcal L_{1,5},\quad Y\cong\mcal L_{1,2},\quad Z\cong\mcal L_{1,4},\quad W\cong\mcal L_{1,3}. \]
They have the first (non-unitary) quantum dimensions. Their conformal dimensions
\[ (h_{1,5},h_{1,2},h_{1,4},h_{1,3})=(11,\frac98,\frac{53}8,\frac{10}3) \]
match our first conformal dimensions mod 1. Thus, our classification result immediately implies that the massive RG flow described by $\mcal B_A$ has $\text{GSD}\in\{4,5\}$ and $\mcal B$ symmetry is spontaneously broken.

\appendix
\setcounter{section}{0}
\renewcommand{\thesection}{\Alph{section}}
\setcounter{equation}{0}
\renewcommand{\theequation}{\Alph{section}.\arabic{equation}}


\begin{thebibliography}{30}
\bibitem{EGNO15}
  P.~Etingof, S.~Gelaki, D.~Nikshych and V.~Ostrik, ``Tensor Categories,'' American Mathematical Society, 2015.
\bibitem{KK23etale}
  K.~Kikuchi,
  ``Classification of connected \'etale algebras in pre-modular fusion categories up to rank three,''
  [arXiv:2311.15631 [math.QA]].
\bibitem{KO01}
  A.~Kirillov Jr. and V.~Ostrik, ``ON A q-ANALOG OF THE MCKAY CORRESPONDENCE AND THE ADE CLASSIFICATION OF slb2 CONFORMAL FIELD,'' Advances in Mathematics 171(2002), 183-227. https://doi.org/10.1006/aima.2002.2072 [arXiv:math/0101219 [math.QA]].
\bibitem{BD11}
  T.~Booker and A.~Davydov, ``Commutative Algebras in Fibonacci Categories,'' Journal of Algebra 355(2012), 176-204. https://doi.org/10.1016/j.jalgebra.2011.12.029
  [arXiv:1103.3537 [math.CT]].
\bibitem{EP09}
   D.~E.~Evans and M.~Pugh, ``$SU(3)$-Goodman-de la Harpe-Jones subfactors and the realization of $SU(3)$ modular invariants,'' Rev. Math. Phys. 21 (2009), 877–928. https://doi.org/10.1142/S0129055X09003761
   [arXiv:0906.4252 [math.OA]].
\bibitem{KL02}
  Y.~Kawahigashi and R.~Longo,
  ``Classification of local conformal nets: Case c \ensuremath{<} 1,''
  Annals Math. \textbf{160}, 493-522 (2004)
  [arXiv:math-ph/0201015 [math-ph]].
\bibitem{EM21}
  C.~Edie-Michell, ``TYPE II QUANTUM SUBGROUPS OF slN . I: SYMMETRIES OF LOCAL MODULES,'' Communications of the American Mathematical Society 3(2023), 112-165. https://doi.org/10.1090/cams/19
  [arXiv:2102.09065 [math.QA]].
\bibitem{G23}
  T.~Gannon,
  ``Exotic quantum subgroups and extensions of affine Lie algebra VOAs -- part I,''
  [arXiv:2301.07287 [math.QA]].
\bibitem{CEM23}
  D.~Copeland and C.~Edie-Michell, ``CELL SYSTEMS FOR Rep(Uq(slN )) MODULE CATEGORIES,''
  [arXiv:2301.13172 [math.QA]].
\bibitem{O00}
   A.~Ocneanu, ``The classification of subgroups of quantum $SU(N)$,'' Quantum Symmetries in Theoretical Physics and Mathematics R. Coquereaux et al (ed), American Mathematical Society, Providence, 2002, pp.133–159.
\bibitem{CZW18}
  S.~X.~Cui, M.~S.~Zini and Z.~Wang, ``On generalized symmetries and structure of modular categories,'' Sci China Math, 62(2019), 417–446. https://doi.org/10.1007/s11425-018-9455-5
  [arXiv:1809.00245 [math.QA].
\bibitem{GK94}
  D.~Gepner and A.~Kapustin,
  ``On the classification of fusion rings,''
  Phys. Lett. B \textbf{349}, 71-75 (1995)
  doi:10.1016/0370-2693(95)00172-H
  [arXiv:hep-th/9410089 [hep-th]].
\bibitem{RSW07}
  E.~Rowell, R.~Stong, Z.~Wang,
  ``On classification of modular tensor categories,'' Communications in Mathematical Physics 292(2009), 343-389. https://doi.org/10.1007/s00220-009-0908-z
  [arXiv:0712.1377 [math.QA]].
\bibitem{BNRW15}
  P.~Bruillard, S.H.~Ng, E.C.~Rowell, and Z.~Wang, ``ON CLASSIFICATION OF MODULAR CATEGORIES BY RANK,'' International Mathematics Research Notices 2016(2016), 7546-7588. https://doi.org/10.1093/imrn/rnw020
  [arxiv:1507.05139 [math.QA]].
\bibitem{DMNO10}
  A.~Davydov, M.~Müger, D.~Nikshych and V.~Ostrik, ``The Witt group of non-degenerate braided fusion categories,'' Journal für die reine und angewandte Mathematik (Crelles Journal), 2013(677), 135-177. https://doi.org/10.1515/crelle.2012.014
  [arXiv:1009.2117 [math.QA]].
\bibitem{P95}
  B.~Pareigis, ``On Braiding and Dyslexia,'' Journal of Algebra 171(1995), 413-425. https://doi.org/10.1006/jabr.1995.1019
\bibitem{O01}
  V.~Ostrik, ``MODULE CATEGORIES, WEAK HOPF ALGEBRAS AND MODULAR INVARIANTS,'' Transformation Groups \textbf8, 177–206 (2003). https://doi.org/10.1007/s00031-003-0515-6
  [arXiv:math/0111139 [math.QA]].
\bibitem{KK23GSD}
  K.~Kikuchi,
  ``Ground state degeneracy and module category,''
  [arXiv:2311.00746 [hep-th]].
\bibitem{ENO02}
  P.~Etingof, D.~Nikshych and V.~Ostrik,
  ``On fusion categories,''
  Annals of Mathematics 162, no. 2 (2005): 581–642. http://www.jstor.org/stable/20159926
  [arXiv:math/0203060 [math.QA]].
\bibitem{O02}
  V.~Ostrik, ``Fusion categories of rank 2,'' Mathematical Research Letters 10 (2002): 177-183. https://dx.doi.org/10.4310/MRL.2003.v10.n2.a5
  [arXiv:math/0203255 [math.QA]].
\bibitem{O13}
  V.~Ostrik, ``Pivotal fusion categories of rank 3,''  Mosc. Math. J., 15(2015), 373–396. https://doi.org/10.17323/1609-4514-2015-15-2-373-396
  [arXiv:1309.4822 [math.QA]].
\bibitem{LPR20}
  Z.~Liu, S.~Palcoux and Y.~Ren, ``Classification of Grothendieck rings of complex fusion categories of multiplicity one up to rank six,'' Lett Math Phys 112, 54 (2022). https://doi.org/10.1007/s11005-022-01542-1
  [arXiv:2010.10264 [math.CT]].
\bibitem{VS22}
  G.~Vercleyen and J.~Slingerland, ``On Low Rank Fusion Rings,''
  [arXiv:2205.15637 [math-ph]].
\bibitem{L14}
  H.~K.~Larson, ``Pseudo-unitary non-self-dual fusion categories of rank 4,'' Journal of Algebra 415(2014), 184-213.
  https://doi.org/10.1016/j.jalgebra.2014.05.032
  [arXiv:1401.1879 [math.QA]].
\bibitem{DZD16}
  J.~C.~Dong, L.~Y.~Zhang and L.~Dai, ``Non-trivially graded self-dual fusion categories of rank 4,'' Acta Mathematica Sinica, English Series 34(2018), 275-287.
  https://doi.org/10.1007/s10114-017-6375-0
  [arXiv:1603.03125 [math.RA]].
\bibitem{anyonwiki}
  ``AnyonWiki,'' https://anyonwiki.github.io/
\bibitem{D90}
  P.~Deligne, ``Catégories tannakiennes,'' The Grothendieck Festschrift, Vol. II. Progr. Math. 87(1990), 111–195.
\bibitem{K96}
  A.~Yu.~Kitaev, ``Quantum error correction with imperfect gates,'' In Proceedings of the Third International Conference on Quantum Communication and Measurement, September 25-30, 1996.
\bibitem{K97}
  A.~Y.~Kitaev,
  ``Fault tolerant quantum computation by anyons,''
  Annals Phys. \textbf{303}, 2-30 (2003)
  doi:10.1016/S0003-4916(02)00018-0
  [arXiv:quant-ph/9707021 [quant-ph]].
\bibitem{TY98}
  D.~Tambara and S.~Yamagami, ``Tensor categories with fusion rules of self-duality for finite abelian groups,'' Journal of Algebra 209(1998), 692-707 https://doi.org/10.1006/jabr.1998.7558
\bibitem{BSS02}
  F.~A.~Bais, B.~J.~Schroers and J.~K.~Slingerland,
  ``Broken quantum symmetry and confinement phases in planar physics,''
  Phys. Rev. Lett. \textbf{89}, 181601 (2002)
  doi:10.1103/PhysRevLett.89.181601
  [arXiv:hep-th/0205117 [hep-th]].
\bibitem{BSS02'}
  F.~A.~Bais, B.~J.~Schroers and J.~K.~Slingerland,
  ``Hopf symmetry breaking and confinement in (2+1)-dimensional gauge theory,''
  JHEP \textbf{05}, 068 (2003)
  doi:10.1088/1126-6708/2003/05/068
  [arXiv:hep-th/0205114 [hep-th]].
\bibitem{BS08}
  F.~A.~Bais and J.~K.~Slingerland,
  ``Condensate induced transitions between topologically ordered phases,''
  Phys. Rev. B \textbf{79} (2009), 045316
  doi:10.1103/PhysRevB.79.045316
  [arXiv:0808.0627 [cond-mat.mes-hall]].
\bibitem{TW19}
  R.~Thorngren and Y.~Wang,
  ``Fusion Category Symmetry I: Anomaly In-Flow and Gapped Phases,''
  [arXiv:1912.02817 [hep-th]].
\bibitem{HLS21}
  T.~C.~Huang, Y.~H.~Lin and S.~Seifnashri,
  ``Construction of two-dimensional topological field theories with non-invertible symmetries,''
  JHEP \textbf{12}, 028 (2021)
  doi:10.1007/JHEP12(2021)028
  [arXiv:2110.02958 [hep-th]].
\bibitem{FMS}
  P.~Di Francesco, P.~Mathieu and D.~Senechal,
  ``Conformal Field Theory,'' doi:10.1007/978-1-4612-2256-9
\bibitem{KK21}
  K.~Kikuchi,
  ``Symmetry enhancement in RCFT,''
  [arXiv:2109.02672 [hep-th]].
\bibitem{KK22free}
  K.~Kikuchi,
  ``Emergent symmetry and free energy,''
  [arXiv:2207.10095 [hep-th]].
\bibitem{Z90}
  A.~B.~Zamolodchikov,
  ``Thermodynamic Bethe Ansatz in Relativistic Models. Scaling Three State Potts and Lee-yang Models,''
  Nucl. Phys. B \textbf{342}, 695-720 (1990)
  doi:10.1016/0550-3213(90)90333-9
\bibitem{Z91}
  A.~B.~Zamolodchikov,
  ``From tricritical Ising to critical Ising by thermodynamic Bethe ansatz,''
  Nucl. Phys. B \textbf{358}, 524-546 (1991)
  doi:10.1016/0550-3213(91)90423-U
\bibitem{M91}
  M.~J.~Martins,
  ``The Thermodynamic Bethe ansatz for deformed W A(N)-1 conformal field theories,''
  Phys. Lett. B \textbf{277}, 301-305 (1992)
  doi:10.1016/0370-2693(92)90750-X
  [arXiv:hep-th/9201032 [hep-th]].
\bibitem{RST94}
  F.~Ravanini, M.~Stanishkov and R.~Tateo,
  ``Integrable perturbations of CFT with complex parameter: The M(3/5) model and its generalizations,''
  Int. J. Mod. Phys. A \textbf{11}, 677-698 (1996)
  doi:10.1142/S0217751X96000304
  [arXiv:hep-th/9411085 [hep-th]].
\bibitem{DDT00}
  P.~Dorey, C.~Dunning and R.~Tateo,
  ``New families of flows between two-dimensional conformal field theories,''
  Nucl. Phys. B \textbf{578}, 699-727 (2000)
  doi:10.1016/S0550-3213(00)00185-1
  [arXiv:hep-th/0001185 [hep-th]].
\end{thebibliography}
\end{document}